\documentclass[11pt,dvips]{article}
\usepackage{latexsym}
\usepackage{srcltx}
\usepackage{amsmath,amsthm,amsfonts,amssymb}
\usepackage[T1]{fontenc}
\usepackage[latin1]{inputenc}
\usepackage{aeguill}
\usepackage{url}
\usepackage{enumerate}
\usepackage{graphicx}
\usepackage[a4paper,textwidth=15cm,textheight=25cm]{geometry}
\usepackage{export}
\usepackage{srcltx}

\newtheorem{thm}{Theorem}[section]

\newtheorem*{thm*}{Theorem}
\newtheorem{dfn}[thm]{Definition} 
\newtheorem*{dfn*}{Definition}
\newtheorem{dfnbis}{Definition}
\newtheorem{cor}[thm]{Corollary}
\newtheorem*{cor*}{Corollary}

\newtheorem{prop}[thm]{Proposition} 
\newtheorem*{prop*}{Proposition} 
\newtheorem*{properties*}{Properties} 
 
\newtheorem{lem}[thm]{Lemma} 
\newtheorem*{lem*}{Lemma}

\newtheorem*{claim*}{Claim} 
 
\newtheorem*{fact*}{Fact}

\newtheorem*{qst*}{Question}

\newtheorem*{pb*}{Problem}

\theoremstyle{remark}
 \newtheorem*{const*}{Construction}

\newtheorem*{rem*}{Remark}
\newtheorem{rem}[thm]{Remark}
\newtheorem*{example*}{Example}

\newtheorem{example}[thm]{Example}

\newenvironment{SauveCompteurs}[1]{%
\newcommand{\monparametre}{#1}
\openexport{\monparametre_sauve}
  \Export{thm}\Export{section}\Export{subsection}\Export{subsubsection}
\closeexport}{}

\newenvironment{UtiliseCompteurs}[1]{%
\newcommand{\monparametre}{#1}
\openexport{\monparametre_aux}
  \Export{thm}\Export{section}\Export{subsection}\Export{subsubsection}
\closeexport
\PackageInfo{export}{\MessageBreak
Importations from \monparametre_sauve.xpt\MessageBreak}%
\InputIfFileExists{\monparametre_sauve.xpt}{\relax}{\relax}%
\renewcommand{\label}[1]{}
}{%
\PackageInfo{export}{\MessageBreak
Importations from \monparametre_aux.xpt\MessageBreak}%
\InputIfFileExists{\monparametre_aux.xpt}{\relax}{\relax}}

\newlength{\espaceavantspecialthm}
\newlength{\espaceapresspecialthm}
{\normalfont \vskip \espaceapresspecialthm}

\newenvironment{specialthm*}[1]{
\vskip\espaceavantspecialthm \noindent \textbf{#1} \itshape}%
{\normalfont \vskip \espaceapresspecialthm}

\newlength{\espaceavantenonce}
\newlength{\espaceapresenonce}
\newcommand{\fontetitreun}[1]{\textbf{#1}} 
\newcommand{\fontetitredeux}[1]{\textit{#1}} 

{\normalfont \vskip \espaceapresenonce}

\newenvironment{enonce1*}[1]{
\vskip\espaceavantenonce \noindent \fontetitreun{#1} \itshape}%
{\normalfont \vskip \espaceapresenonce}

{\vskip \espaceapresenonce}

\newenvironment{enonce2*}[1]{
\vskip\espaceavantenonce \noindent \fontetitredeux{#1} }%
{\vskip \espaceapresenonce}

\makeatletter
\edef\@tempa#1#2{\def#1{\mathaccent\string"\noexpand\accentclass@#2 }}
\@tempa\rond{017}
\makeatother

\newcommand{\es}{\emptyset}
\renewcommand{\phi}{\varphi} 
\newcommand{\m} {^{-1}}

\newcommand {\ra} {\rightarrow}

\newcommand{\actson}{\curvearrowright}

\newcommand{\dunion}{\sqcup}

\newcommand{\Dunion}{\bigsqcup} 

\newcommand{\ie} {i.e.\ }

\newcommand {\cala} {{\mathcal {A}}}   
\newcommand {\calb} {{\mathcal {B}}}   
   
\newcommand {\cald} {{\mathcal {D}}}

\newcommand {\calh} {{\mathcal {H}}}

\newcommand {\calk} {{\mathcal {K}}}   
   
\newcommand {\calm} {{\mathcal {M}}}

\newcommand {\calp} {{\mathcal {P}}}

\newcommand {\cals} {{\mathcal {S}}}   
\newcommand {\calt} {{\mathcal {T}}}   
\newcommand {\calu} {{\mathcal {U}}}

\newcommand {\calz} {{\mathcal {Z}}}

\newcommand {\bbH} {{\mathbb {H}}}

\newcommand {\bbZ} {{\mathbb {Z}}}   

\newcommand{\Out} {\mathop{\mathrm{Out}}}
\newcommand{\Aut} {\mathop{\mathrm{Aut}}}

\newcommand {\Z} {{\mathbb {Z}}}

\newcommand{\inc}{\subset}
\newcommand{\bo}{\partial}

\newcommand\elli{{\mathrm{ell}}}

\usepackage{pdfsync}

\begin{document}

\title{JSJ decompositions: definitions, existence, uniqueness. \\
I: The JSJ deformation space.}
\author{Vincent Guirardel, Gilbert Levitt}
\date{}

\maketitle

\begin{abstract}
We give a simple general definition of JSJ decompositions by means of a universal maximality property.
The JSJ decomposition should not be viewed as a tree (which is not uniquely defined) but as a canonical deformation space of trees. We prove that JSJ decompositions of  finitely presented groups always exist,
 without any   assumption  on edge groups.
Many examples are given.
\\

\textbf{This paper and its companion \url{arXiv:1002.4564} have been replaced by \url{arXiv:1602.05139}.} 
 \end{abstract}

\section{Introduction }

JSJ decompositions first appeared in 3-dimensional topology with the theory of the characteristic submanifold
by Jaco-Shalen and Johannson \cite{JaSh_JSJ,Johannson_JSJ}.
For an orientable irreducible   closed $3$-manifold $M$, this can be described as follows.
Let  $\calt\subset M$ be a maximal disjoint union of non-parallel embedded tori
such that any immersed torus can be homotoped to be disjoint from them.
Then $\calt$ is unique up to isotopy, and any connected component of $M\setminus\calt$ is either atoroidal,
or a Seifert fibered manifold.

This was carried over  to group theory by Kropholler \cite{Kro_JSJ}
for some Poincar\'e duality groups of dimension at least $3$, and by Sela for torsion-free hyperbolic groups \cite{Sela_structure}.
The initials \emph{JSJ}, standing for  {\bf J}aco-{\bf S}halen and {\bf J}ohannson, were popularized by Sela, 
and constructions of JSJ decompositions were  given in more general settings by many authors 
\cite{RiSe_JSJ,Bo_cut,DuSa_JSJ,FuPa_JSJ,DuSw_algebraic,ScSw_regular}.

In this group-theoretical context, one has a group $G$ and a class of subgroups $\cala$ (such as cyclic groups, abelian groups...), 
and one tries to understand splittings   (i.e.\   graph of groups decompositions)  of $G$ over groups in $\cala$
(in  \cite{DuSw_algebraic,ScSw_regular}, one looks at almost invariant sets rather than splittings, 
in closer analogy to the $3$-manifold situation).
The family of tori $\calt$ of the $3$-manifold is replaced by a splitting  of $G$ over $\cala$.
The   authors construct   a canonical splitting 
 enjoying a long list of properties, rather  specific to each case.

These ideas 
have had a vast influence and range of applications, 
from the isomorphism problem and the structure of the group of automorphisms of hyperbolic groups, 
to diophantine geometry over groups  \cite{Sela_structure,Sela_isomorphism,Sela_diophantine1}.

In this paper, we do not propose a construction of JSJ decompositions,
but rather  a simple abstract \emph{definition}  stated by means of a  universal maximality property,  together with general 
  \emph{existence} and \emph{uniqueness} statements stated in terms of \emph{deformation spaces} (see below).
 
 The JSJ decompositions constructed in \cite{RiSe_JSJ,Bo_cut,DuSa_JSJ,FuPa_JSJ} are JSJ decompositions in our sense,
and we view these constructions as \emph{descriptions} of JSJ decompositions and  their flexible vertices.
The regular neighbourhood  of \cite{ScSw_regular} is of a different nature, whose relation with usual JSJ decompositions is explored in \cite{GL5}.

Several results of this paper and its sequel \cite{GL3b} were announced in \cite{GL_barcelona}.

\paragraph{A universal property.}

To explain the definition, let us  first consider free decompositions of a group $G$,
\ie decompositions of $G$ as the fundamental group of  a graph of groups with trivial edge groups, 
or equivalently   actions of $G$ on a tree $T$ with trivial edge stabilizers.

If $G=G_1*\dots *G_n$ where each $G_i$ is non-trivial, non-cyclic,  and freely indecomposable, the Bass-Serre tree $T_0$ of this decomposition is \emph{maximal}: if $T$ is  associated to any other free decomposition, then $T_0$ \emph{dominates} $T$ in the sense that there is a $G$-equivariant map $T_0\to T$. In other words, 
among free decompositions of $G$,  the tree 
$T_0$ is as  far  as possible from the trivial  tree (a point):
its vertex stabilizers are as small as possible (they are conjugates of the $G_i$'s). The maximality condition does not determine $T_0$ uniquely, we will come back to this key fact later.

\begin{figure}[htbp]
   \centering
   \includegraphics{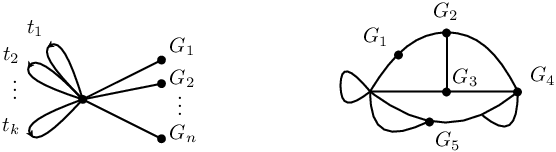}
   \caption{Graph of groups decompositions corresponding to two Grushko trees}
   \label{fig_grushko}
 \end{figure}

 If now  $G=G_1*\dots * G_p *F_q$ with $G_i$ as above and $F_q$  a free group, one can take $T_0$ to be the universal covering of   one of the 
graphs of groups
 pictured on Figure \ref{fig_grushko}.  Its vertex stabilizers are precisely the conjugates of the $G_i$'s (if $G$ is free, its action on $T_0$ is free). We call such a tree a Grushko tree.

When more general decompositions are allowed, there may not exist a maximal tree. The  fundamental example is the following. 
Consider a closed surface $\Sigma$, and  two simple closed curves $c_1,c_2$ in $\Sigma$ with non-zero intersection number. 
Let $T_i$ be the Bass-Serre tree of the 
 associated 
splitting of $\pi_1(\Sigma)$ over $\bbZ\simeq\pi_1(c_i)$.
Since $c_1$ and $c_2$ have positive intersection number, $\pi_1(c_1)$ is hyperbolic in $T_2$ (it does  not fix a point) and vice-versa.
Using one-endedness of  $\pi_1(\Sigma)$, it is an easy exercise to check that there is no splitting of $\pi_1(\Sigma)$
which dominates both $T_1$ and $T_2$. In this case 
  there is no hope of having a maximal splitting over cyclic groups, similar to $T_0$ in the case of free splittings.

To overcome this difficulty, one restricts to  \emph{universally elliptic} splittings, defined  as follows.
Let $\cala$ be a class of subgroups of $G$, stable under taking subgroups and under conjugation.
We only consider  $G$-trees with edge stabilizers in $\cala$, which we call  \emph{$\cala$-trees}.

\begin{dfnbis}
 An   $\cala$-tree is \emph{universally elliptic} if its edge stabilizers are elliptic in every $\cala$-tree.
\end{dfnbis}

Recall that $H$ is elliptic in $T$ if it fixes a point in $T$. Free decompositions are universally elliptic, but the trees $T_1,T_2$ introduced above are not.

If an   $\cala$-tree $T$ is universally elliptic, then one can \emph{read} any $\cala$-tree with finitely generated edge stabilizers
from $T$ by blowing up   vertices, and then performing a finite sequence of collapses and folds on this blow-up (see Remark \ref{rem_fold}).

\begin{dfnbis}
  A \emph{JSJ decomposition} (or JSJ tree) of $G$ over $\cala$ is an $\cala$-tree $T$ such that 
  \begin{itemize}
  \item $T$ is universally elliptic;
  \item $T$ dominates any other universally elliptic tree.
  \end{itemize}
\end{dfnbis}

The second condition is a maximality condition. It means that, if $T'$ is universally elliptic, there is an equivariant map $T\to T'$. Equivalently, vertex stabilizers of $T$ are elliptic in every universally elliptic tree.

If $\cala$ consists of all subgroups with a given property (being cyclic, abelian, slender, ...), we refer to, say, cyclic trees, cyclic JSJ decompositions when working over $\cala$.

\paragraph{Uniqueness.}
JSJ trees are not unique. Returning to the example of free decompositions, one obtains trees 
with the same maximality property as $T_0$ by precomposing the action of $G$ on $T_0$ with any automorphism of $G$. One may also change the topology of the quotient graph $T_0/G$ (see Figure \ref{fig_grushko}). The canonical object is   not a single tree, but the set of all trees with trivial edge stabilizers and non-trivial vertex stabilizers conjugate to the $G_i$'s, a \emph{deformation space}.

\begin{dfnbis}[Forester]
The \emph{deformation space} of a 
tree   $T$ is the set of 
trees 
$T'$ such that $T$ dominates $T'$ and $T'$ dominates $T$. Equivalently, two 
trees are in the same deformation space $\cald$ if and only if they have the same elliptic subgroups. 
\end{dfnbis}

More generally,  given a family of subgroups   $\widetilde\cala\inc\cala$, 
one considers deformation spaces over $\widetilde\cala$ by restricting to trees in $\cald$ with edge stabilizers in $\widetilde\cala$.

For instance, Culler-Vogtmann's outer space (the set of free actions of $F_n$ on trees) is a deformation space. Just like outer space, any deformation space   is a complex in a natural way, and it is contractible (see \cite{GL2,Clay_contractibility}).

If $T$ is a JSJ tree, a  tree $T'$ is a JSJ tree if and only  if  $T'$ is universally elliptic,  $T$ dominates $T'$, and $T'$ dominates $T$. 
In other words, $T'$ should belong to the   deformation space of $T$ over $\cala_{ell}$, where $\cala_{ell}$ is the family of universally elliptic groups in $\cala$.

\begin{dfnbis}
If non-empty, the set of all JSJ decompositions of $G$ is a deformation space over $\cala_{ell}$ called the \emph{JSJ deformation space} (of $G$ over $\cala$). We denote it by $\cald_{JSJ}$.
\end{dfnbis}

The canonical object is therefore not a particular JSJ decomposition, but the   JSJ deformation space.

It is a general fact that two trees belong to the same deformation space if and only if one can pass from one to the other by applying a finite sequence of moves of certain types (see \cite{For_deformation,GL2,For_splittings,ClFo_Whitehead}  and Remark \ref{moves}). 
The statements about uniqueness of the JSJ up to certain moves which appear in \cite{RiSe_JSJ,DuSa_JSJ,FuPa_JSJ} are special cases of this general fact.

If $\cala$ is invariant under the group of automorphisms of $G$ (in particular if $\cala$ is defined by restricting the isomorphism type), then so is $\cald_{JSJ}$, and one can gain information about $\Aut(G)$ and $\Out(G)$ by studying their action on the contractible complex $\cald_{JSJ}$ \cite{CuVo_moduli,McCulloughMiller_symmetric,GL1,Clay_deformation}.

There are nice situations  (for instance, splittings of a one-ended hyperbolic group over two-ended subgroups \cite{Bo_cut}) 
when one can construct a preferred   tree in $\cald_{JSJ}$ (or in a related deformation space). Such a tree is a fixed point for the action of  $\Out(G)$. 
In \cite{GL4,GL3b}, we   explain how   fixed points may sometimes be constructed as trees of cylinders or as maximal decompositions encoding compatibility of splittings (see the discussion at the end of this introduction).    

\paragraph{Existence.}
Existence of the JSJ deformation space over any class $\cala$ of subgroups (without any smallness assumption) when $G$ is finitely presented
is a simple consequence of Dunwoody's accessibility.

\begin{UtiliseCompteurs}{thm_exist}  
\begin{thm}
If $G$ is finitely presented, 
  the JSJ deformation space $\cald_{JSJ}$ of $G$ over $\cala$ exists. It contains a  tree whose edge and vertex stabilizers are finitely generated.
\end{thm}
\end{UtiliseCompteurs}

In the sequel to this paper \cite{GL3b} we use
  acylindrical accessibility to construct and describe  JSJ decompositions  of certain finitely generated groups, for instance abelian decompositions of CSA groups.

\paragraph{Description.}
Given $G$ and $\cala$, all JSJ trees have the same vertex stabilizers, provided we restrict to stabilizers not in $\cala$. 
A vertex $v$, or its stabilizer $G_v$, is 
  \emph{rigid} if $G_v$ is universally elliptic (\ie it fixes a point in every $\cala$-tree). 
For instance, all vertices of the  Grushko tree $T_0$ studied above are rigid. 
  
But the essential feature of JSJ theory is the description of \emph{flexible} vertices (those which are not rigid), in particular  the fact that 
  in many contexts (\cite{RiSe_JSJ,DuSa_JSJ,FuPa_JSJ}; see Theorem \ref{thm_lesgens})
 flexible vertex stabilizers are extensions of $2$-orbifold groups with boundary,
attached to the rest of the group in a particular way (quadratically hanging subgroups). 
In other words, the example of   trees $T_1,T_2$ given above using intersecting curves on  a surface is  often 
the only source of flexibility.

 We prove in \cite{GL3b} that this is also the case for flexible subgroups of the JSJ deformation space
of a  relatively hyperbolic group with small parabolic subgroups over the class of its small subgroups.
 
However,   in certain natural 
 situations, flexible groups are not   quadratically hanging subgroups, 
for instance for JSJ decompositions over abelian groups  (see Subsection \ref{ab}).

  We also consider \emph{relative} JSJ decompositions, where one imposes that finitely many finitely generated  subgroups be elliptic in all splittings considered, and we show that the  description of flexible subgroups remains valid in this context (Theorem \ref{thm_lesgensrel}).

\begin{pb*}
  Describe flexible vertices of the JSJ deformation space of a finitely presented group over small subgroups.\\
\end{pb*}

In a sequel to this paper  \cite{GL3b}, we 
produce  a canonical tree instead of a  canonical   deformation space, by 
replacing universal ellipticity  by the more rigid notion of universal compatiblity.
If $T$ is universally elliptic, and $T'$ is arbitrary, one can refine $T$ (by blowing up vertices) to a tree $\hat T$  which dominates $T'$. 
We say that $T$ is \emph{universally compatible} if, given $T'$, one can refine $T$ to a tree $\hat T$ which \emph{refines} $T'$. 
We show that, if $G$ is finitely presented, there exists a maximal deformation space containing a universally compatible tree. 
Unlike $\cald_{JSJ}$, this deformation space always contains a preferred tree $T_{co}$, the \emph{compatibility JSJ tree.}
This is somewhat similar to Scott and Swarup's  construction who
construct a canonical tree $T_{SS}$ which is compatible with (\ie which \emph{encloses}) any almost invariant set \cite{ScSw_regular}.  The tree $T_{co}$ dominates $T_{SS}$, sometimes strictly.

Additionnally, we use acylindrical accessibility and trees of cylinders
to produce and describe JSJ deformation spaces and compatibility
JSJ trees for some classes of \emph{finitely generated} groups including CSA groups, $\Gamma$-limit groups with $\Gamma$ a hyperbolic group  (possibly with torsion),
and  relatively hyperbolic groups. This also includes relative JSJ decompositions,  relative to an arbitrary  (possibly infinite) family of subgroups.
\\

Let us now describe the contents of this paper. After preliminary sections, the JSJ deformation space is defined and constructed in Section \ref{exi}.
We also explain there why the constructions of \cite{RiSe_JSJ,DuSa_JSJ,FuPa_JSJ} are JSJ decompositions in our sense. Section \ref{sec_relative} is devoted to relative JSJ decompositions. In Section \ref{sec_rig}, we give examples of JSJ decompositions with all vertices rigid  
(like   the Grushko trees $T_0$ in the discussion above, these trees dominate every tree  under consideration). 
This includes splittings over finite groups, parabolic splittings of  relatively hyperbolic groups (relative to parabolic subgroups), 
and locally finite trees with small edge  stabilizers (in particular, generalized Baumslag-Solitar groups).

In Section \ref{sec_QH}, we define quadratically hanging (QH) subgroups, and give some general properties of   these groups.
We state results from \cite{RiSe_JSJ,DuSa_JSJ,FuPa_JSJ} 
showing that flexible subgroups are QH-subgroups in JSJ decompositions over cyclic groups, or two-ended groups, or slender groups (Theorem \ref{thm_lesgens}).
We also show that, under suitable assumptions, any QH vertex stabilizer of any $\cala$-tree is elliptic in JSJ trees (Section \ref{sec_QH_ell}). 

Finally, in Section \ref{sec_flex}, 
we  give examples of flexible subgroups over abelian groups which are not QH-subgroups. 
We also show how a relative JSJ decomposition of a group $G$ may be viewed as an absolute decomposition of  
a larger group $\hat G$, using a filling construction. 
This allows us to extend the description of flexible subgroups  to the relative case  (the original proofs probably extend to the
relative case, but this is unlikely to ever appear in print). 
\\

\setcounter{tocdepth}{2}
\tableofcontents
\section{Preliminaries}  \label{prel}

In all of this paper, $G$ will be a finitely generated group. Sometimes finite presentation will be needed, for instance
to prove existence of JSJ decompositions.

We also  fix a family $\cala$ of subgroups of $G$ which is  stable under
conjugation and under taking subgroups. If $H$ is a subgroup, we sometimes write $\cala_{ | H}$ for the set of subgroups of $H$ which belong to $\cala$. 

\subsection*{Trees}
A \emph{splitting} of $G$ is an isomorphism of $G$ with the fundamental group of a graph of groups.
A one-edge splitting (a graph of groups with one edge) is a splitting as an amalgam or an HNN extension.
Using Bass-Serre theory, we view a splitting as an action of $G$ on a tree $T$ without inversion.
The action $G\actson T$ is \emph{trivial} if $G$ fixes a point, and \emph{minimal} if there is no proper $G$-invariant subtree.

In this paper, all trees 
are endowed with a minimal action of $G$ without inversion (we allow the trivial case when  $T$ is a point).
We identify two  trees if there is an equivariant isomorphism between them. 

An
\emph{$\cala$-tree} is a tree $T$ whose edge stabilizers belong to $\cala$. We sometimes say
that
$T $ is
\emph{over $\cala$}. We also say cyclic tree (abelian tree, ...) when $\cala$ is the family of cyclic (abelian, ...) subgroups (we consider the trivial group as cyclic).

An element  or a subgroup  of $G$ is \emph{elliptic} in $T$ if it fixes a point.
 If $H_1\inc H_2\inc G$ with $H_1$ of finite index in $H_2$, then $H_1$ is elliptic if and only if $H_2$ is. 
An element which is not elliptic is \emph{hyperbolic}, it  has a unique axis on which it acts by translation.
If a finitely generated subgroup of $G$   is not elliptic, it contains a hyperbolic element and has a unique minimal invariant subtree.  

Let $H$ be a subgroup of $G$ containing no hyperbolic element. If it does not fix a point, then it is not  finitely generated,  and it fixes a unique end of $T$: there is a ray $\rho$ such that each finitely generated subgroup of $H$ fixes a subray of $\rho$.

A  tree $T$ is \emph{irreducible} if there exist  two hyperbolic elements $g,h\in G$ whose axes intersect in a compact set.
Equivalently, $T$ is not irreducible if and only if $G$ fixes a point, or an end, or preserves a line of $T$.

We denote by $V(T)$ and $E(T)$ the set of vertices and (non-oriented) edges of $T$ respectively, by  
 $G_v$ or $G_e$ the stabilizer of a vertex $v$ or an edge $e$.  

If $v\in V(T)$, its \emph{incident edge groups} are the stabilizers of edges containing $v$, viewed as subgroups of $G_v$. We denote by $\calp_v$ the set of incident edge groups. It is a finite union of conjugacy classes of subgroups of $G_v$.
\subsection*{Relative trees}

 Besides $\cala$, we sometimes also   fix an arbitrary   set $\calh$
 of subgroups of $G$, and we restrict to  
$\cala
$-trees   
$T$
 such that each $H\in\calh$ is elliptic in $T$  (in terms of graphs of groups, $H$ is contained in a conjugate of a vertex group). We call such a tree an \emph{$(\cala, \calh)$-tree}, or a tree over $\cala$ relative to $\calh$. 
 The set of $(\cala, \calh)$-trees does not change if we enlarge $\calh $ by making it invariant under conjugation. 
 
 In particular, suppose that $G_v$ is a vertex stabilizer of a tree. 
We often consider splittings of $G_v$ relative to $\calp_v$, the set of incident edge groups. 
 Such a splitting extends to a splitting of $G$ (see Remark   \ref{rem_bas}).
 
\subsection*{Morphisms, collapse maps, refinements}

Maps between trees will always be $G$-equivariant and linear on edges.
 We  mention  two particular classes of maps.

A map $f:T\ra T'$ between two trees is a \emph{morphism} if
each edge of $T$ can be written as a  finite union of
subsegments, each of which is mapped
bijectively onto a segment in $T'$. Equivalently, $f$ is a
morphism if and only if one may subdivide
$T$ and $T'$ so that $f$ maps each edge onto an edge (in particular, no edge of $T$ is collapsed). 
Folds are examples of morphisms (see \cite{BF_bounding}). 

A \emph{collapse} map $f:  T\to T'$ is a map obtained by collapsing   edges to points (by equivariance, the set of collapsed edges is $G$-invariant). 
Equivalently, $f$ preserves alignment: the image of any arc $[a,b]$ is a point or  the arc $[f(a),f(b)]$.
We say that $T'$ is a  collapse of $T$, and that $T$ is a \emph{refinement} of $T'$. 
 In terms of graphs of groups, one obtains $T'/G$ by collapsing edges in $T/G$.

\subsection*{Domination and deformation spaces}\begin{dfn}
A  tree $T_1$ \emph{dominates} $T_2$ if there is an equivariant map from $T_1$ to $T_2$.
Equivalently, $T_1$ dominates $T_2$ if every vertex stabilizer of $T_1$ fixes a point in $T_2$: every subgroup which is elliptic in $T_1$ is also elliptic in $T_2$.  
\end{dfn}

In particular, every refinement of $T_1$ dominates $T_1$.  Beware that domination is defined by considering ellipticity of subgroups, not just of elements.

\begin{dfn}\label{dfn_deformation}
Having the same elliptic subgroups is an equivalence relation on the set of $\cala$-trees, an equivalence class is called  a \emph{deformation space} over $\cala$.
\end{dfn}

In other words, if $T_0$ is an $\cala$-tree, the deformation space of $T_0$ is the set of all $\cala$-trees $T$ such
that $T$ and $T_0$ dominate each other, or equivalently the set of  $\cala$-trees having the same elliptic subgroups as $T_0$.
We say that a deformation space $\cald$ dominates a space $\cald'$ if trees in $\cald$ dominate those of $\cald'$. 
The deformation space of the trivial tree is called the \emph{trivial deformation space.} It is the only deformation space in which $G$ is elliptic. 

\begin{rem} \label{moves} The deformation space of $T_0$ is the space of trees that may be obtained from $T_0$ by applying a finite sequence of collapse and expansion moves 
(\cite{For_deformation}, see also \cite{ClFo_Whitehead}). When the deformation space $\cald$ is non-ascending as defined in \cite{GL2}, 
for instance when all groups in $\cala$ are finite, any two reduced trees in $\cald$ may be joined by a finite sequence of slide moves. 
These facts will not be used in the present paper.  
\end{rem}

\subsection*{Slender, small subgroups}

A group  $H$  is \emph{slender} if $H$ and all its subgroups are finitely generated.
Examples of slender groups include virtually polycyclic groups.
We denote by $VPC_n$ (resp.\ $VPC_{\leq n}$) the class of virtually polycyclic groups
of Hirsch length exactly  (resp.\ at most) $n$. If a slender group acts on a tree, it fixes a point or leaves a line invariant.

Following  \cite[page 177]{GL2}, we say that 
a subgroup $H\subset G$ (possibly infinitely generated) is \emph{small in $\cala$-trees} 
if, 
whenever  $G$ acts on an $\cala$-tree, $H$ fixes a point, or an end,  or leaves a line invariant.
Any small subgroup (in the sense of \cite{BF_bounding}), in particular any subgroup not containing $F_2$, is small in $\cala$-trees. 
Unlike smallness, smallness in 
$\cala$-trees is stable under  taking subgroups. 
It is also a commensurability invariant  (but it is not invariant under abstract isomorphism). 
  Moreover, $H$ is small in $\cala$-trees if and only if all its finitely generated subgroups are.

\subsection*{Accessibility}

\renewcommand{\labelenumi}{(\theenumi)}
 
Constructions of JSJ decompositions are based on accessibility theorems stating that, given  suitable $G$ and $\cala$, there is 
an a priori 
bound for the number of orbits of edges of $\cala$-trees,  
under the assumption that there is no redundant vertex: if $v$ has valence 2, it is the unique fixed point of some $g\in G$.
This holds in particular:
\begin{enumerate} \item if $G$ is finitely presented and all groups in $\cala $ are finite \cite{Dun_accessibility};
\item if $G$ is finitely presented, all groups  in $\cala $ are small, and the trees are  \emph{reduced} in the sense of \cite{BF_bounding};
\item if $G$ is finitely generated and all groups in $\cala $ are finite with bounded order \cite{Linnell};
\item if $G$ is finitely generated and the trees are $k$-acylindrical for some $k$ \cite{Sela_acylindrical};
\item if $G$ is finitely generated and the trees are $(k,C)$-acylindrical: the stabilizer of any segment of length $>k$ has order $\le C$  \cite{Weidmann_accessibility} (\cite{Delzant_accessibilite} for finitely presented groups). 
\end{enumerate}
 
In this paper, we use a version of Dunwoody's accessibility given in \cite{FuPa_JSJ} (see Proposition \ref {prop_accessibility}).
 In \cite{GL3b} we use acylindrical accessibilty.

\section{Universal ellipticity}

Fix $\cala$ as above. Let $T_1,T_2$ be two $\cala$-trees.

\begin{dfn}[Ellipticity of trees]
 $T_1$ is \emph{elliptic} with respect to $T_2$ if every edge stabilizer
of $T_1$ is elliptic in $T_2$.
\end{dfn}

When $T_1$ is  {elliptic} with respect to $T_2$, one   can \emph{read $T_2$ from $T_1$}:
  there is a refinement
$\Hat T_1$ of $T_1$ which dominates $T_2$. More precisely:

\begin{lem}\label{lem_refinement} 
Let $T_1,T_2$ be minimal  $\cala$-trees.
If  $T_1$ is elliptic with respect to $T_2$, there is a minimal  $\cala$-tree   $\Hat T_1$    such that:
\begin{enumerate}
\item $\Hat T_1$ is a refinement of $T_1$ and dominates $T_2$;
\item the stabilizer of any edge of $\Hat T_1$ fixes an edge in $T_1$ or in $T_2$;
\item every edge stabilizer of $T_2$ contains an edge stabilizer of $\Hat T_1$;
\item a subgroup of $G$ is elliptic in $\Hat T_1$ if and only if it is elliptic in both $T_1$ and
$T_2$.
\end{enumerate}
\end{lem}

\begin{rem} Conversely, if  some refinement   $\Hat T_1$ dominates $T_2$, then $T_1$ is  {elliptic} with respect to $T_2$: edge stabilizers of $T_1$ are elliptic in $\Hat T_1$,
hence also   in $T_2$.
\end{rem}

\begin{rem}\label{rem_fold}
   If edge stabilizers of $T_2$ are finitely generated, then $T_2$ can be obtained from
$\Hat T_1$ by a finite number of collapses and folds \cite{BF_bounding}.
\end{rem}

\begin{proof}
For each vertex $v\in V(T_1)$, choose   any
$G_v$-invariant subtree
$Y_v$ of $T_2$ (possibly $ T_2$ itself).
For each   edge $e=vw\in E(T_1)$, 
choose     vertices $p_v\in Y_{v}$ and $p_w\in Y_w$ fixed by $G_e$; this is possible
because
$G_e$ is elliptic in
$T_2$.   We make these choices $G$-equivariantly.

We can now define a  tree $\Hat T_1$
by blowing up each vertex $v$ of $T_1$ into $Y_v$, and attaching     edges   of $T_1$ using the
points $p_v$. Formally, we consider the disjoint union $ \Dunion_{v\in V(T_1)}Y_v$, and for
each edge $ vw$ of
$T_1$ we attach an edge between $p_v\in Y_v$ and $p_w\in Y_w$.

Clearly, $\Hat T_1$ is a refinement of $T_1$  and dominates $T_2$: one recovers
$ T_1$ from $\Hat T_1$    by collapsing each $Y_v$ to a point, and one constructs a map  
$f:\Hat T_1\to T_2$ by mapping each added edge  $vw$ to the segment $[p_v,p_w]$. If $\Hat T_1$ is not minimal, we replace it by its minimal subtree. We now consider stabilizers.

The stabilizer
of any edge  of
$\Hat T_1$ fixes an edge in $T_2$ or in $T_1$, depending on whether the edge  lies in some $Y_v$ or not.  In particular, $\Hat T_1$ is an $\cala$-tree. By minimality, any edge $e$ of $T_2$ is contained in the image of some edge of $\Hat T_1$, so $G_e$ contains an edge stabilizer of $\Hat T_1$.

If a subgroup $H\subset G$ is elliptic in $\Hat T_1$, then it is elliptic in $T_1$ and $T_2$.
Conversely, if $H$ is elliptic in $T_1$ and $T_2$,
let $v\in T_1$ be fixed by $H$. Since $H$ fixes a point in $Y_v$,
it is elliptic in $\Hat T_1$.
\end{proof}

\begin{rem}\label{rem_bas}
The basic idea of  this construction may be summarized as follows. If $v$ is a vertex of a tree $T$, and $T_v$ is a tree with an action of $G_v$ such that for every incident edge $e=vw$  of $T$ the group $G_e$ is elliptic in $T_v$, one can refine $T$ by replacing $v$ by $T_v$ (and performing a similar replacement at every vertex in the orbit of $v$).
\end{rem}

The deformation space of   $\Hat T_1$ only depends on the  deformation spaces $\cald_1$ and $\cald_2$  of $T_1$ and $T_2$: it   is the
lowest deformation space dominating $\cald_1$ and $\cald_2$. If $T_1$ dominates $T_2$ (resp.\ $T_2$ dominates $T_1$),
then
$\hat T_1$ is in the same deformation space as $T_1$ (resp.\ $T_2$).
Moreover,  there is some
symmetry: if $T_2$ also happens to be elliptic with respect to $T_1$,  then $\Hat T_2$ is in the
same deformation space as $\Hat T_1$.

\begin{lem} \label{cor_Zor}
\begin{enumerate}
\item  If $T_1$ refines $T_2$ and does not belong to the same deformation space, some $g\in G$ is hyperbolic in $T_1$ and elliptic in $T_2$.
\item
If  $T_1$ is elliptic with respect to  $T_2$, and every $g\in G$ which is elliptic in $T_1$ is also elliptic in $T_2$, then $T_1$ dominates $T_2$.  
\end{enumerate}
\end{lem}

\begin{proof} 
One needs only   prove the first assertion when $T_2$ is obtained from $T_1$ by collapsing the orbit of an edge $e=uv$.
If $u$ and $v$ are in the same orbit, or if $G_e\neq G_u$ and $G_e\neq G_v$, then some hyperbolic element of $T_1$
becomes elliptic in $T_2$. Otherwise, $T_1$ and $T_2$ are in the same deformation space.

For the second assertion, assume that $T_1$ does not dominate $T_2$. Then the  tree  $\Hat T_1$ given by
Lemma
\ref{lem_refinement} does not belong to the same deformation space as $T_1$. 
Since it is a refinement of $T_1$, some $g\in G$ is elliptic in $T_1$ and hyperbolic in $\Hat T_1$. By Assertion  (4) of Lemma
\ref{lem_refinement}, this element is hyperbolic in $T_2$, a contradiction.  Thus $T_1$ dominates $T_2$.\end{proof}

\begin{dfn}[Universally elliptic]
A subgroup $H\subset G$ is \emph{universally elliptic} (over $\cala$) if 
it is elliptic in every  $\cala$-tree.
An $\cala$-tree $T$ is \emph{universally elliptic} if its edge stabilizers are universally elliptic, \ie if $T$ is elliptic with respect to every $\cala$-tree.
\end{dfn}

\begin{lem}\label{lem_sup}
\begin{enumerate}
\item   If $T_1,T_2 $ are universally elliptic,   the refinement  $\Hat T_1$ given by
Lemma
\ref{lem_refinement} is   universally elliptic. 
\item
If there is a morphism $f:S\to T$, and $T $ is
universally elliptic, then so is $S$  (see Section \ref{prel} for the definition of morphisms).  
\item If $T_1,T_2$ are universally elliptic and have the same elliptic elements, they belong to the same deformation space.
\end{enumerate}
\end{lem}

\begin{proof}   The first assertion follows from Assertion   (2)  of Lemma \ref{lem_refinement}.
The second assertion is clear. The third one follows from the second assertion of Lemma \ref{cor_Zor}.
\end{proof}

The following lemma will be useful in \cite{GL3b}. 
\begin{lem} \label{lem_Zorn} 
Let $(T_i )_{i \in I}$ be any family of   trees. There exists a countable subset
$J\subset I$ such that, if 
$T$ is elliptic  with respect to  every $T_i$ ($i\in I$),
and $T$ dominates every $T_{j}$ for $j \in J$, 
then $T$ dominates $T_i $ for all $i\in I$.
\end{lem}

\begin{proof} Since $G$ is countable, we can   find a countable  $J$ such that, if an element $g\in G$ is
hyperbolic in some $T_i$, then it is hyperbolic in some $T_{j}$ with $j \in J$. 
If $T$ dominates every $T_{j }$ for $j \in J$, any $g$ which is elliptic in $T$ is elliptic in every $T_i$. 
By Lemma  \ref{cor_Zor}, the tree $T$ dominates every $T_i$. 
\end{proof}

For many purposes, it is enough to consider one-edge splittings, \ie $\cala$-trees with only one orbit of edges.

\begin{lem}\label{lem_oneed}
 Let $S$ be an $\cala$-tree. \begin{enumerate}
\item
$S$ is universally elliptic if and only if
it is elliptic with respect to every one-edge splitting (over $\cala$). 
\item
$S$ dominates every universally elliptic  $\cala$-tree if and only if it dominates every universally elliptic one-edge
splitting.
\end{enumerate}
\end{lem}

\begin{proof} By induction on the  number of orbits of
edges, using the following lemma.
\end{proof}

\begin{lem}
Let $T$ be a  tree, and $H$ a subgroup of $G$.
  Let $E_1\dunion E_2$ be a partition of $E(T)$ into two $G$-invariant sets.
Let $T_1,T_2$ be the trees obtained from $T$ by collapsing $E_1$ and $E_2$ respectively.
\begin{enumerate}
\item If a subgroup $H$ is elliptic in $T_1$ and $T_2$, then $H$ is elliptic in $T$.
\item If a tree $T'$ dominates $T_1$ and $T_2$, then it dominates $T$.
\end{enumerate}
\end{lem}

\begin{proof}
Let $x_1\in T_1$ be a vertex fixed by $H$. Let $Y\subset T$ be its preimage. It is a subtree
because the map $T\ra T_1$ is a collapse map.
Now $Y$ is $H$-invariant and embeds into $T_2$.
Since $H$ is elliptic in $T_2$, it fixes a point in $Y$, therefore it is elliptic in $T$.  One shows (2) 
by applying (1) to vertex stabilizers of $T'$.
\end{proof}

\section{The JSJ deformation space}  \label{exi}

\subsection{Definitions}

We fix $\cala$ as above. 
  We define $\cala_\elli\subset \cala$ as the set of   groups in $\cala$ which are
universally elliptic ($\cala_\elli$ is stable under conjugating and taking subgroups). An $\cala$-tree is universally elliptic if and only if it is an $\cala_\elli$-tree.

\begin{dfn}[JSJ deformation space]\label{prop_def}
A deformation space $\cald_{JSJ}$ of $\cala_\elli$-trees which is maximal for  
domination is called \emph{the JSJ deformation space} of $G$ over $\cala$ (it is unique if it exists by the first assertion of Lemma \ref{lem_sup}).

Trees in $\cald_{JSJ}$ are called \emph{JSJ trees}, or  \emph{JSJ decompositions}, of $G$ over $\cala$. 
They are precisely those $\cala$-trees $T $ which are universally
elliptic, and   which dominate every
 universally elliptic tree. 
\end{dfn}


   In general there are
many JSJ trees, but they  all belong to    the same deformation space
and therefore have a lot in common (see Section 4 of \cite{GL2}). In
particular \cite[Corollary 4.4]{GL2},  they have the same vertex stabilizers, except possibly  for 
vertex stabilizers in $\cala_{\elli}$  (the groups in $\cala$ which are universally elliptic). 

\begin{dfn} [Rigid and flexible vertices] Let $H=G_v$ be a vertex stabilizer   of  a JSJ  tree $T$ over $\cala$.  
We say that $H$ is \emph{rigid} if it is universally elliptic,
\emph{flexible} if it is not. 
We also say that the vertex $v$ is rigid  (flexible). If $H$ is flexible, we say that it  is a  flexible  subgroup of $G$
over $\cala$.  
 \end{dfn} 

The definition of  flexible subgroups of $G$ does not
depend on the choice of the  JSJ tree $T$. 
The heart of JSJ  theory is to understand   flexible groups.   They will be discussed in
Sections \ref{sec_QH} and \ref{sec_flex}.

\subsection{Existence}

\begin{SauveCompteurs}{thm_exist}
  \begin{thm}\label{thm_exist_mou}
    Assume that $G$ is finitely presented.  Then the JSJ deformation
    space $\cald_{JSJ}$ of $G$ over $\cala$ exists. It contains a tree
    whose edge and vertex stabilizers are finitely generated.
  \end{thm}
\end{SauveCompteurs}

There is no   hypothesis, such as smallness or finite generation, on
the elements of $\cala$.
Finite presentability of $G$ is used to prove the existence of $\cald_{JSJ}$, through the
following version of   Dunwoody's accessibility  whose proof will be  given after that of Theorem \ref{thm_exist_mou}.

\begin{prop}[Dunwoody's accessibility]\label{prop_accessibility}
  Let $G$ be finitely presented.
  Assume that $T_1\leftarrow\dots\leftarrow T_k \leftarrow T_{k+1}\leftarrow \dots$
is a sequence of refinements of   $\cala$-trees.
There exists an $\cala$-tree $S $ such that:
\begin{enumerate}
\item for   $k$ large enough, there is a morphism $S\ra T_k$ (in particular, $S$ dominates $T_k$);
\item  each edge and vertex  stabilizer of $S$ is finitely generated.
\end{enumerate}
\end{prop}

\begin{rem*}
Note that the maps $T_{k+1}\ra T_k$ are required to be \emph{collapse maps}.
\end{rem*}

One may view the following standard result as a special case (apply Proposition \ref
{prop_accessibility} to the constant sequence $T\leftarrow T\leftarrow\dots$).

\begin{cor}\label{cor_fg}
  If $G$ is finitely presented, and   $T$ is an $\cala$-tree,
  there exists  a morphism $f:S\ra T$ where $S$ is an $\cala$-tree   with finitely generated
edge and  vertex stabilizers (if $T$ is universally elliptic, so is $S$ by Lemma  \ref{lem_sup}). \qed
\end{cor}

 \begin{proof}[Proof of Theorem \ref{thm_exist_mou}]
Let $\calu $ be the set of minimal universally elliptic $\cala$-trees 
with finitely generated edge and vertex stabilizers, up to equivariant isomorphism.
It is non-empty since it contains the trivial tree. An element of $\calu$ is described by a
finite graph of groups with finitely generated edge and vertex groups. Since  $G$  only has
countably many finitely generated subgroups, and there are countably many 
 homomorphisms from a given finitely generated  group to another, the set $\calu$ is
countable.  

  By Corollary \ref{cor_fg}, it suffices  to produce a universally elliptic $\cala$-tree dominating every $U\in\calu$.
Choose an enumeration $\calu=\{U_1,U_2,\dots,U_k,\dots\}$.
 We define inductively a universally elliptic $\cala$-tree $T_k$ which refines $T_{k-1}$ and
dominates $U_1,\dots, U_k$ 
 (it may have infinitely generated edge or vertex stabilizers).
We start with $T_1=U_1$. Given $T_{k-1}$ which dominates $U_1,\dots,U_{k-1}$,
we let $T_k$ be a   universally elliptic refinement of $T_{k-1}$ which dominates $U_k$,
given by Lemma  \ref{lem_sup}. Then $T_k$ also dominates $U_1,\dots, U_{k-1}$
because
$T_{k-1}$ does.

  Apply  Proposition \ref{prop_accessibility} to the sequence $T_k$. The tree $S$ is universally
elliptic  because there are morphisms $S\to T_k$.
Furthermore $S$ dominates every $T_k$, hence every $U_k$. 
This shows that $S$ is a JSJ tree over $\cala$.
\end{proof}

\begin{proof}[Proof of Proposition \ref{prop_accessibility}]

This is basically Proposition 5.12 of \cite{FuPa_JSJ}. We sketch the argument  for completeness. 

Let  $X$ be a   Cayley  simplicial 2-complex for $G$, 
\ie $X$ is a simply connected simplicial 2-complex   with a free
cocompact action of
$G$. 
By considering preimages of midpoints of edges under
  a suitable equivariant map $X\to T_k$, one constructs a $G$-invariant pattern $P_k$ on $X$ such that there is a morphism  $S_{P_k}\ra T_k$, where  $S_{P_k}$ is the tree   dual to $P_k$.  
Since $T_{k+1}$ refines $T_k$, one can assume    $P_{k }\subset P_{k+1}$. 

By Dunwoody's theorem \cite[Theorem 2.2]{Dun_accessibility},  there is a bound on the number of non-parallel tracks in $X/G$.
Thus  there exists $k_0$ such that for all $k\geq k_0$ one obtains $S_{P_k}$  from
$S:=S_{P_{k_0}}$ by subdividing   edges, 
and therefore there is a morphism $S
\ra T_k$ for $k\geq k_0$. 
Edge and vertex stabilizers of $S $ are finitely generated since they are images in
$G$ of fundamental groups of nice compact subsets  of $X/G$.  
\end{proof}

\begin{rem} \label{rem_linn}
When $\cala$ consists of finite groups with bounded order, Proposition \ref{prop_accessibility} is true if $G$ is only assumed to be finitely generated. This follows from Linnell's accessibility \cite{Linnell}: for $k$ large, $T_{k+1}$ is just a subdivision of  $T_k$. The JSJ deformation space therefore always exists in this case. 
\end{rem}

We record the following simple facts for future reference. 
\begin{lem} \label{lem_rigid}
Any $\cala$-tree $T$ with
universally elliptic \emph{vertex} stabilizers is a JSJ tree.
\end{lem}

\begin{proof} This is clear since $T$ dominates every $\cala$-tree.
\end{proof}

\begin{lem}\label{lem_rafin}
  If  there is  a JSJ tree, then any $\cala$-tree $T$ with
universally elliptic edge stabilizers may be refined  to a JSJ tree.
\end{lem}

\begin{proof} 
Let $T_2$ be a JSJ tree. Apply Lemma \ref{lem_refinement}
with $T_1=T$. The refinement  $\hat T_1$ of $T$ is universally elliptic and dominates $T_2$,
so is a JSJ tree.
\end{proof}

\subsection{Relation with other constructions}
\label{sec_gens}
Several authors have constructed JSJ splittings of finitely presented groups in various
settings. We explain here why those splittings are JSJ splittings in our sense (results in the literature are often
stated only for one-edge splittings, but this is not a restriction by  Lemma \ref{lem_oneed}).

In \cite{RiSe_JSJ}, Rips and Sela consider cyclic splittings of a one-ended group $G$ (so
$\cala$ consists of all  cyclic subgroups of $G$,  including the trivial group). Theorem 7.1 in
\cite{RiSe_JSJ} says that the produced JSJ splitting  is universally elliptic (this is statement
(iv)) and maximal (statement (iii)).  The uniqueness up to deformation is statement (v). 

In \cite{FuPa_JSJ}, Fujiwara and Papasoglu consider all   splittings of a   group
over the class $\cala$ of its slender subgroups.
Statement (2) in \cite[Theorem 5.13]{FuPa_JSJ}  says that the produced JSJ splitting is elliptic
with respect to any minimal splitting. By Proposition 3.7 in \cite{FuPa_JSJ}, any
splitting is dominated  by a minimal splitting, so universal ellipticity holds.
Statement (1) of Theorem 5.15   in \cite{FuPa_JSJ}  implies maximality. 

In  the work of Dunwoody-Sageev  \cite{DuSa_JSJ}, the authors consider splittings of a  
group $G$ over slender subgroups in a class $\calz\calk$ 
such that $G$ does not split over finite
extensions of infinite index subgroups of $\calz\calk$ (there are restrictions on the class
$\calz\calk$,  but one can typically take $\calz\calk=VPC_n$,
see \cite{DuSa_JSJ} for details). 
In our notation, $\cala$ is the set of subgroups of elements of $\calz\calk$.
Universal ellipticity of the constructed splitting follows from statement (3) in the Main Theorem of  \cite{DuSa_JSJ}, and from the fact that 
any edge group is contained in a white vertex group.
Maximality follows from the fact that white vertex groups are universally elliptic (statement (3)) 
and that black vertex groups   either are in $\calz\calk$
(in which case they are universally elliptic by the non-splitting assumption made on $G$), or are
$\calk$-by-orbifold groups  and hence are necessarily elliptic in any JSJ tree (see Proposition \ref{qhe} below).

\section{Relative JSJ decompositions}\label{sec_relative}

\subsection{Definition and existence}\label{rela}

 The construction of the JSJ deformation space may be done in a relative setting. In this subsection, we show how to adapt the  
 construction given in Section \ref{exi}. 
  In Section \ref{sec_flex},  we will give another approach (valid under similar hypotheses) and study flexible vertices. 
In \cite{GL3b} we will construct relative JSJ decompositions  of certain groups
using acylindrical accessibility instead of Dunwoody's accessibility (this does not require   finite presentability).

 Besides $\cala$, we   fix a   set $\calh$
 of subgroups of $G$ and we only
consider  $(\cala, \calh)$-trees: 
$\cala
$-trees   
$T$
 such that each $H\in\calh$ is elliptic in $T$.    Of
course, universal ellipticity is defined using only $(\cala, \calh)$-trees. We  now denote by  $\cala_\elli$ the
set of elements of $\cala$ which are elliptic in every  $(\cala, \calh)$-tree. 

The \emph{JSJ deformation space} of $G$ over $\cala$ relative to 
$\calh$, if it exists, is  the  unique   deformation space of $(\cala_\elli, \calh)$-trees which is maximal for  
domination.  Note that the JSJ space does not change if we enlarge $\calh $ by making it invariant under conjugation. 

 Suppose that $\calh=\{H_1,\dots,H_p\}$ is a finite family.
One says that $G$  is  \emph{finitely presented relative to $\calh$} if it  is  a quotient 
 $G\simeq (F *\tilde H_1*\dots * \tilde H_p)/N$, where $F$ is a finitely generated free group, $N$ is the normal
closure of a finite set, and each $\tilde H_i$   maps isomorphically onto
$H_i$.
If we drop the finiteness condition on $N$, we say that $G$ is \emph{finitely generated relative
to $\calh$}. 
  Note that replacing each $H_i$ by a conjugate does not change the fact that $G$ is relatively finitely presented or generated.
A finitely presented group is finitely presented relative to any finite collection of finitely generated subgroups.  
A relatively hyperbolic group is finitely presented relative to its maximal  parabolic subgroups \cite{Osin_relatively}.

\begin{thm}\label{thm_exist_mou_rel} 
Assume that  $G$ is finitely presented relative to $\calh=\{H_1,\dots,H_p\}$.
Then the JSJ deformation space $\cald_{JSJ}$ of $G$ over $\cala$ relative to 
$\calh$ exists. 

Moreover, $\cald_{JSJ}$ contains a tree with finitely generated edge stabilizers
(and finitely generated vertex stabilizers if all $H_i$'s are finitely generated).
\end{thm}

The theorem is proved as in the non-relative case, using    a relative version of
Proposition  \ref{prop_accessibility}: the trees $T_k$ and $S$ are now $(\cala , \calh)$-trees, and 
stabilizers of $S$ are  finitely generated relative to finitely many conjugates of elements of $\calh$ (hence finitely generated if all groups in $\calh$ are finitely generated).

The proof of this relative accessibility is similar to the proof of Proposition \ref{prop_accessibility},
except that one has  to work with a \emph{relative Cayley 2-complex} $X$ as follows. 
Consider connected simplicial complexes $X_1,\dots, X_p$ (not necessarily compact),  with base points $*_i$, such that   $\pi_1(X_i,*_i)=H_i$. Glue finitely many $2$-cells 
to  the wedge of the $X_i$'s and finitely many circles  to get a space $X$ with
$\pi_1(X)=G$. In the universal cover  $\pi:\Tilde X\ra X$, consider $Y=\pi\m(X_1)\dunion \dots \dunion \pi\m(X_p)$.
The stabilizer of each connected component of $Y$ is conjugate to some $H_i$, and $\pi(\Tilde X\setminus Y)$ is relatively compact.
 If each $H_i$ fixes a point in a tree $T$, there is an equivariant map $f:\Tilde X\ra T$ that maps each component of
$Y$ to a vertex of $T$, and preimages of midpoints of edges of $T$ define a $G$-invariant pattern in $\Tilde X\setminus Y$.
The tree dual to this pattern is relative to $\calh$ because the pattern does not intersect $Y$.
One concludes as in the absolute case using Dunwoody's bound on the number of non-parallel tracks in $X\setminus (X_1\cup\dots \cup X_p)$.

\subsection{JSJ decompositions of vertex groups}\label{sec_JSJ_vx}
We  now record   relations between the JSJ decomposition  of a group, and relative splittings of vertex stabilizers.
All trees are assumed to have   edge stabilizers in $\cala$.
Splittings of a vertex stabilizer $G_v$ are over $\cala_{ | G_v}$, the family of subgroups of $G_v$ belonging to $\cala$,  and relative to the incident edge stabilizers.  

\begin{lem}\label{lem_flex}
  Let $ G_v$ be a vertex stabilizer   of  a universally elliptic  tree $T$. Let $\calp_v$ be the set of incident edge stabilizers. 
\begin{enumerate}
\item A subgroup of $ G_v$ is universally elliptic (as a subgroup of $G$) if and only if it is elliptic in every splitting of $G_v$ relative to $\calp_v$. 
\item Assume that $T$ is a JSJ tree. Then $G_v$ has no  non-trivial splitting  relative to $\calp_v$ over a universally elliptic subgroup.
 The group $G_v$ is universally elliptic  if and only if it has no     non-trivial  splitting  relative to $\calp_v$.
\item Let $S$ be a JSJ tree of $G$. The action of $G_v$ on its minimal subtree in $S$ is a JSJ tree of $ G_v$   relative to $\calp_v$.
 In particular, $G_v$ is elliptic in the JSJ deformation space   of $G$ if and only if its  JSJ decomposition relative to $\calp_v$   is trivial.
\end{enumerate}
\end{lem}

\begin{proof} 
If  $H\inc G_v$ is not   elliptic in some tree   with an action of $G$, restricting the action to $G_v$ shows the ``if'' direction of the first assertion; note that the groups of $\calp_v$ are elliptic in the induced splitting because they are universally elliptic. By Remark \ref{rem_bas}, any splitting of $G_v$ relative to $\calp_v$ may be  used  to refine $T$.  This shows the other direction.

 The second assertion  follows similarly, using 
 the maximality of JSJ trees.  For the  third  assertion, we may assume by Lemma   \ref{lem_rafin} that $S$ is a refinement of $T$. We then use Assertion  (1).
 \end{proof}

\begin{lem}\label{lem_JSJrel} 
Let $T$ be a universally elliptic tree.
Then $G$ has a JSJ tree if and only if  every vertex stabilizer 
$G_v $  of $T$ has a JSJ tree relative to $\calp_v$.

In this case, one can refine  $T$ using these  trees so as to get  a JSJ tree of $G$.
\end{lem}
 
\begin{proof}
If $G$ has a JSJ tree, so does $G_v$ by Lemma \ref{lem_flex} (3).
Conversely, assume that $G_v $ has a JSJ tree $T_v$  
relative to $\calp_v$.
Choose such a $T_v$ in a $G$-equivariant way.
Let $\Hat T$ be the corresponding refinement of $T$ as in Remark \ref{rem_bas}.
It is universally elliptic by Assertion  (1) of Lemma \ref{lem_flex}.
Maximality of $\Hat T$ follows from maximality of $T_v$  as in the proof of  Lemma \ref{lem_rafin}: 
given a universally elliptic $T_2$, the refinement of $\Hat T$ which dominates $T_2$ belongs to the same deformation space as $\Hat T$ by  Assertion  (1) of Lemma \ref{lem_flex}.
\end{proof}

\begin{rem}\label{rem_flexrel}
Lemmas \ref{lem_flex} and  \ref{lem_JSJrel}  hold  if  one considers   JSJ trees relative to $\calh$, as in Subsection \ref{rela}, 
provided one adds to $\calp_v$ all subgroups of $G_v$   of the form $H\cap G_v$ with $H$ conjugate to a group of $\calh$. Since $H\cap G_v$ fixes an edge incident to $v$ if $H$ is not contained in $G_v$, it suffices to add to $\calp_v$ the subgroups of $G_v$  which are  conjugate to a group of $\calh$.
\end{rem}

\section{Rigid examples} \label{sec_rig}

In this section,    we shall construct $\cala$-trees whose  vertex
stabilizers  are  universally elliptic; 
they are JSJ trees by  Lemma \ref {lem_rigid} (they dominate every $\cala$-tree), and all their vertices are rigid. 
Unless specified otherwise,  we only assume that $G$ is finitely generated, so  JSJ trees are not guaranteed to exist.

\subsection{Free groups}

Let $G=F_n$ be a finitely generated free group,
and let $\cala$ be arbitrary.
Then the JSJ deformation space of $F_n$ over $\cala$ is the space of free actions (unprojectivized Culler-Vogtmann's 
Outer Space \cite{CuVo_moduli}).

More generally, if $G$ is virtually free and 
$\cala$ contains all finite subgroups, then $\cald_{JSJ}$ is the space of trees with finite vertex
stabilizers.

\subsection{Free splittings: the Grushko deformation space} \label{freep}
Let $\cala$ consist only of the trivial subgroup
of $G$. Then the JSJ deformation space exists, it
is the outer space introduced in \cite{GL1} (see \cite{McCulloughMiller_symmetric} when no free factor of $G$ is $\bbZ$, \cite{CuVo_moduli} when $G=F_n$).
This is the set of Grushko  trees $T$ in the following sense:
edge stabilizers are trivial and vertex stabilizers are freely indecomposable (we consider 
$\bbZ$ as freely decomposable since it splits over the trivial group). 
 Denoting by 
$G=G_1*\dots*G_p*F_q $   a  decomposition of $G$ given by Grushko's theorem (with $G_i$ non-trivial and
freely indecomposable, and $F_q$ free), the quotient graph of groups $T/G$ 
 is homotopy equivalent to a wedge of $q$ circles, 
 it has one vertex  
 with group $G_i$ for each $i$,  and all other vertex groups are trivial  (see Figure \ref{fig_grushko}).

\subsection{Splittings over finite groups:  the Stallings-Dunwoody  deformation space}
\label{Dunw}
If $\cala$ is the set
of finite subgroups of $G$, 
the JSJ deformation space is the set of  trees whose edge groups are finite and whose vertex
groups have 0 or 1 end 
(by Stallings's theorem,
an $\cala$-tree is maximal for domination if and only if its vertex stabilizers have at most one end).
By definition, the 
JSJ deformation space exists if and only if $G$ is accessible, in particular if $G$ is finitely
presented by Dunwoody's original accessibility result
\cite{Dun_accessibility}.

As mentioned in Remark  \ref{rem_linn},   Linnell's accessibility implies that 
the JSJ space exists if $G$ is finitely generated and $\cala$ consists of finite groups with bounded order.

\subsection{Parabolic splittings}
\label{parab}

Assume that all groups in $\cala$ are universally elliptic. Then all vertex stabilizers of a JSJ tree $T$ are rigid by Lemma \ref{lem_flex} (and Remark \ref{rem_flexrel} in the relative case).

In particular, assume that $G$ is hyperbolic relative to a family of finitely generated subgroups $\calh=\{H_1,\dots,H_p\}$.   Recall  that a subgroup of $G$ is
\emph{parabolic} if it is contained in a conjugate of an $H_i$. We let $\calp$ be the family of parabolic subgroups, and we consider splittings over $\calp$ relative to $\calh$ (equivalently, relative to $\calp$).

Parabolic subgroups are universally elliptic (relative to $\calh$!), so
the  JSJ trees over parabolic subgroups, relative to $\calh$, do not have flexible vertices (these trees exist  
because $G$ is finitely presented relative to $\calh$, see Subsection \ref{rela}).

\subsection{Splittings of small  groups}\label{sec_G_small}

  Recall that $G$ is small in $\cala$-trees if it has no irreducible action on a $\cala$-tree (irreducible meaning non-trivial,
without any invariant line or end).

 \begin{prop} \label{sma}
 If $G$ is  small in  $\cala$-trees, there is at most one non-trivial deformation space containing a universally elliptic $\cala$-tree.
\end{prop}
 
 \begin{proof} If  $T$ is  a non-trivial universally elliptic tree, every vertex stabilizer   contains an edge stabilizer with index at most 2 (the index is 2 if $G$ acts dihedrally on  a line), so is universally elliptic. It follows that any two such trees belong to the same deformation space. 
\end{proof}

If there is a deformation space as in the proposition, it is the JSJ space. Otherwise, the JSJ space is trivial. 

Consider for instance  cyclic splittings of  solvable Baumslag-Solitar groups $BS(1,n)$. If $n=1$ (so $G=\Z^2$), 
there are infinitely many deformation spaces and no non-trivial universally elliptic tree. 
If $n=-1$ (Klein bottle group), there are   exactly two deformation spaces but no non-trivial universally elliptic tree.
If $n\ne\pm1$, the JSJ space is non-trivial, as we shall now see.

 \subsection{Generalized Baumslag-Solitar groups}\label{gbs}

Let $G$ be a generalized Baumslag-Solitar group,
\ie it acts on  a tree $T$ 
with all  vertex and edge stabilizers infinite cyclic.
Let $\cala$ be the set of cyclic subgroups of $G$ (including the trivial subgroup).
Unless $G$ is   isomorphic to  $\Z$, $\Z^2$, or to the   Klein bottle group,
  the deformation space of $T$  is  the JSJ deformation space \cite{For_uniqueness}.

Here is a short proof (the arguments are contained in \cite{For_uniqueness}). We show
  that every vertex stabilizer $H$ of $T$ is universally elliptic.
Clearly, the commensurator of $H$ is $G$. If $H$ is hyperbolic in an $\cala$-tree $T'$, 
its commensurator $G$ preserves its axis,
so $T'$ is a line. This implies that
$G$ is  $\Z$, $\Z^2$, or  a Klein bottle group, a contradiction.

\subsection{Locally finite trees}

We generalize the previous example to small splittings (\ie splittings over small groups). 

We suppose that 
$G$ acts irreducibly on a  locally finite tree $T$ with all edge stabilizers  small in  $\cala$-trees (local finiteness
is equivalent to edge stabilizers having finite index in neighboring vertex stabilizers; in particular,
vertex stabilizers are small in $\cala$-trees). In \cite[Lemma 8.5]{GL2}, we proved that all such trees $T$ belong
to the same deformation space. This happens to be the JSJ
deformation space.

\begin{lem}  Suppose that all groups of $\cala$ are   small in $\cala$-trees.  Then any locally finite irreducible
$\cala$-tree $T$ belongs to the JSJ deformation space over $\cala$.
\end{lem}

\begin{proof} We show that every vertex stabilizer $H$ of  $T$ is universally elliptic.  Assume that $H$
is not elliptic in some $\cala$-tree $T' $.  If $H$ contains a hyperbolic element, then
 it preserves a unique line or end of $T'$ because it is small in $\cala$-trees.  
If it consists of elliptic elements but does not fix a point, then $H$ fixes a unique end of $T'$.
Moreover, in both cases, any finite index subgroup of $H$ preserves the same unique line or end.

As in the previous subsection, local finiteness implies that 
  $G$ commensurizes $H$, so it preserves this $H$-invariant  line or end of $T'$ (in particular, $T'$ is not
irreducible).
We now define a normal subgroup $G'\inc G$ which is small in $\cala$-trees.
If $G$ does not act dihedrally on  a line, we let $G'=[G,G]$.
It  is small in $\cala$-trees: any finitely generated subgroup pointwise fixes a ray of
$T'$, so is contained in an edge stabilizer $G_e\in\cala$.  If $G$ acts dihedrally, we let $G'$ be the kernel of the action.

 Consider the action of the normal subgroup $G'$ on $T$. 
If it is elliptic,   its fixed point set  is $G$-invariant, so by minimality
the action of $G$ factors through the action of an abelian  or dihedral group; this contradicts the irreducibility of
$T$.  If  $G'$ preserves a unique line or end, this line or end is
$G$-invariant because
$G'$ is normal, again    contradicting the irreducibility of $T$.
\end{proof}

\section{QH subgroups}\label{sec_QH}

 Flexible subgroups of the JSJ deformation space are most important, as understanding 
 their splittings conditions the understanding of the splittings of $G$.
In well-understood cases,    flexible subgroups are 
\emph{quadratically hanging (QH) subgroups}.
In this section, after preliminaries about $2$-orbifold groups, we define QH-subgroups     and 
we establish some general properties.
  We then quote several results saying that flexible vertices are QH, and we prove that, under suitable hypotheses,
any QH-subgroup is elliptic in the JSJ deformation space.

 \subsection{2-orbifolds}  \label{orb}

 We consider a compact 2-orbifold $\Sigma$ with $\pi_1(\Sigma)$ not virtually abelian. Such an orbifold is hyperbolic, we may view it as the quotient of a compact orientable hyperbolic surface $ \Sigma_0$ with   geodesic boundary by a finite group of isometries $\Lambda$. If we forget the orbifold structure, $\Sigma$ is homeomorphic to a surface   $\Sigma_{top}$. 
 
 The image of $\bo\Sigma_0$ is the \emph{boundary} $\bo\Sigma$ of $\Sigma$. Each component $C$ of $\bo\Sigma$ is either a component of $\bo\Sigma _{top}$ (a circle) or an arc contained in $\bo\Sigma _{top}$.
The  (orbifold) fundamental group of $C$ is $\Z$ or an infinite dihedral group  accordingly.  
A   \emph{boundary subgroup} is a subgroup of $\pi_1(\Sigma)$ 
which is conjugate to  the fundamental group of a component $C$ of $\partial \Sigma $.

The  closure of the complement of $\bo\Sigma$ in $\bo\Sigma_{top}$ is a union of \emph{mirrors}:
a mirror is the image of a component of the fixed point set of an orientation-reversing element of  $\Lambda$.
Each mirror is itself a circle or  an arc  contained in $\bo\Sigma _{top}$.
Mirrors may be adjacent, whereas boundary components of $\Sigma$ are disjoint.

\begin{lem} \label{arc}
Let $C$ be a boundary component of a compact hyperbolic 2-orbifold $\Sigma$. 
There exists a non-trivial splitting of $\pi_1(\Sigma)$ over $\{1\}$ or $\Z/2\Z$ 
relative to the fundamental groups $J_k$  of all boundary components distinct from $C$.
\end{lem}

\begin{proof}
Any arc $\gamma$ properly embedded in $\Sigma_{top}$ and with endpoints on $C$ defines a free splitting of $\pi_1(\Sigma)$  relative to the groups  $J_k$. In most cases one can choose $\gamma$ so that this splitting is non-trivial. The exceptional cases are when $\Sigma_{top}$ is   a disc or an annulus, and $\Sigma$ has no conical point. 
 
If $\Sigma_{top}$ is   a disc, its boundary circle consists of components of $\bo\Sigma$ and mirrors. 
 Since $\Sigma$ is hyperbolic, there must be a mirror $M$ not adjacent to $C$ 
(otherwise $\bo\Sigma_{top}$ would consist of $C$ and one or two mirrors, or two boundary components and two mirrors). 
An arc $\gamma$ with one endpoint on $C$ and the other on $M$ defines a  splitting over $\Z/2\Z$, which is non-trivial because $M$ is not adjacent to $C$.

If $\Sigma_{top}$ is   an annulus, there are two cases. If  $C$ is an arc, one can find an arc $\gamma$ from $C$ to $C$ as in the general case. If $C$ is a  circle in $\bo\Sigma_{top}$, the other circle contains a mirror  $M$ (otherwise $\Sigma$ would be  a regular annulus) and an arc $\gamma$ from $C$ to $M$ yields a splitting over $\Z/2\Z$.
\end{proof}

\begin{rem}\label{arcb}
 If the splitting constructed is over $K=\Z/2\Z$, this $K$ is contained in a 2-ended subgroup (generated by $K$ and a conjugate).
\end{rem}

We denote by $\tilde \Sigma$ the universal covering of $\Sigma$, a convex subset of $\bbH^2$  with geodesic boundary.  
A  (bi-infinite) geodesic  $\gamma \inc\tilde\Sigma$ is 
  \emph{closed} if its image in $\Sigma$ is compact,
and
\emph{simple} 
if  $h\gamma$ and $\gamma$ are equal or disjoint for all $h\in \pi_1(\Sigma)$. 
If $\gamma\not\subset\bo\tilde\Sigma$, we say that its projection $\delta$ is an \emph{essential simple closed  geodesic} 
in $\Sigma$ (possibly one-sided). We then
denote by $H_\gamma$ the 2-ended subgroup of $\pi_1(\Sigma) $ consisting of elements which preserve $\gamma$ and each of the   half-spaces bounded by    $ \gamma$. 

There is a  non-trivial one-edge splitting of $\pi_1(\Sigma)$ over $  H_\gamma$ 
relative to the    boundary subgroups, we say that it is \emph{dual to $\delta$} 
(if $\delta$ is one-sided, this splitting can be viewed 
as the splitting 
dual to the boundary of a regular neighbourhood of $\delta$, 
a connected $2$-sided simple $1$-suborbifold). 
More generally, any family of disjoint simple closed geodesics $\delta_i$ gives rise to a dual splitting.

\subsection{Quadratically hanging   subgroups}  \label{quah}

Let $Q$ be a vertex stabilizer  of an $\cala$-tree.
 
\begin{dfn}\label{dfn_qh} 
  We say that $Q$ is a \emph{QH-subgroup} (over $\cala$) if it is an
  extension $1\to F\to Q \to \pi_1(\Sigma )\to 1$, with $\Sigma$ a
  hyperbolic 2-orbifold as above, and each incident edge group is an
  \emph{extended boundary subgroup}: its image in $\pi_1(\Sigma)$ is
  finite or contained in a boundary subgroup of $\pi_1(\Sigma)$ (in
  particular, the preimage in $Q$ of the stabilizer of a cone point of
  $\Sigma$ is an extended boundary subgroup).  
  We call $F$ the
  \emph{fiber}.   

 We say that a boundary component $C$ of
  $\Sigma $ is \emph{used} if there exists an incident edge group
  whose image in $\pi_1(\Sigma)$   is contained in $\pi_1(C)$  with  finite index
  (recall that $\pi_1(C)$ is cyclic or dihedral, so having finite
  index is the same as being infinite).
\end{dfn}

Note that the terminology QH is used by Rips-Sela in
  \cite{RiSe_JSJ} with a more restrictive meaning ($F=\{1\}$ and
  $\Sigma$ has no mirror).

\subsection{General properties of QH subgroups}

 Let $\Sigma$ be as above. We first recall that  splittings of $\pi_1(\Sigma)$ over small
  (\ie virtually cyclic) subgroups are dual to families of simple closed curves.

\begin{lem} \label{dual}  Let $T$ be a tree
with a non-trivial minimal action of $\pi_1(\Sigma)$,  without inversion,
  with small edge stabilizers and 
with all boundary subgroups elliptic. 
Then $T$ is  equivariantly isomorphic to the  
tree dual to a family of disjoint 
simple closed geodesics  of $\Sigma$ (possibly one-sided, see Subsection \ref{orb}).

 If  edge stabilizers are not assumed to be small,
$T$ is still dominated by a  
tree dual to a family of disjoint 
simple closed geodesics  of $\Sigma$. 
\end{lem}

\begin{proof} 
When $\Sigma$ is an orientable surface, this follows from 
Theorem III.2.6 of \cite {MS_valuationsI}, noting that all  small subgroups of $\pi_1(\Sigma )$ are cyclic.
If $\Sigma$ is a  $2$-orbifold, we consider a covering  surface $\Sigma_0$ as in Subsection \ref{orb}. 
The action of $\pi_1(\Sigma_0)$ on $T$ is dual to a family of  
closed geodesics on $\Sigma_0$. 
This family is $\Lambda$-invariant and projects to the required family on $\Sigma$. 
The action of $\pi_1(\Sigma)$ on $T$ is dual to this family, as defined in   Subsection \ref{orb}.

The second statement follows from standard arguments (see the proof of   \cite[Theorem III.2.6]{MS_valuationsI}).
\end{proof}

The following proposition shows that Lemma \ref{dual} applies under natural conditions.

\begin{prop} \label{prop_fue}
Let   $Q$ be a QH vertex group  of an $\cala$-tree $T$, with fiber $F$   and   underlying orbifold  $\Sigma$.
\begin{enumerate}
\item \label{it_fue}
If $\cala$ consists of slender groups, and $F$ is slender, 
then $F$ is  universally elliptic.
\item \label{it_factor} If  $F$ is universally elliptic, and if $T'$ is any $\cala$-tree,
 there is a $Q$-invariant subtree $T_Q\inc T'$ such that 
the action of $Q$ on $T_Q$ factors through an action of $\pi_1(\Sigma)$. 
\item \label{it_used}  Assume that $F$, and any subgroup containing $F$ with index 2, belongs to $\cala$. If $T$ is a JSJ decomposition over $\cala$, and if   $F$ is universally elliptic,
then every boundary component of $\Sigma$ is used.  Moreover, every extended boundary subgroup of $Q$ is universally elliptic. 
\end{enumerate}
\end{prop}

\begin{proof}
Suppose $F$ is not   universally
elliptic. Being slender, it would act non-trivially on  a line   in some
$\cala$-tree  with an action of $G$. Since
$F$ is normal in
$Q$, this line would be $Q$-invariant and $Q$ would act on it with slender edge
stabilizers. This is a contradiction since $Q$ is not slender. 

For the second assertion, let $T_Q$ be the fixed subtree of $F$. Since $F$ is normal in $Q$, it is $Q$-invariant.

Let $C$ be a boundary component of $\Sigma$. Lemma \ref{arc} yields 
a non-trivial splitting of $Q$ over a group  
 containing F with index	$\le2$, 
 hence in $\cala$  and universally elliptic. 
If $C$ is not used, every incident edge group is elliptic in this splitting. 
  By Remark \ref{rem_bas}, one can use this decomposition of $Q$ to produce a 
universally elliptic splitting of $G$ in which $Q$ is not elliptic.
This contradicts the maximality of the JSJ decomposition. 
Since $F$ and all edge stabilizers of $T$ are universally elliptic, so are extended boundary subgroups. 
\end{proof}

We now prove that, under natural hypotheses, the only universally elliptic elements in a QH-subgroup $Q$
lie in extended boundary subgroups.

\begin{prop}\label{prop_relper}
  Let   $Q$ be a QH vertex group  of an $\cala$-tree $T$, with fiber $F$   and   underlying orbifold  $\Sigma$.
\begin{enumerate}
\item 
If  $F$, and extended boundary subgroups of $Q$, are universally elliptic, 
but   $Q$ is not universally elliptic, then    $\Sigma$ contains a essential simple closed geodesic (as defined in Subsection \ref{orb}).
\item 
Assume that the preimages in $Q$ of all two-ended subgroups of $\pi_1(\Sigma)$   belong to  $\cala$. 
If $\Sigma$ contains an essential simple closed geodesic, any universally elliptic element of $Q$ lies in an extended boundary subgroup. 
 In particular, $Q$ is not universally elliptic. 
\end{enumerate}
\end{prop}

\begin{proof}
Since $Q$ is not universally elliptic, it  acts non-trivially on some $\cala$-tree $T'$. As in Proposition \ref{prop_fue},
the  action of $Q$ on its minimal subtree 
$T_Q$ factors through $\pi_1(\Sigma)$.
By assumption, the boundary subgroups of $\pi_1(\Sigma)$ fix a point in $T'$,
so the action is dominated by  one which is dual to a system of simple closed essential geodesics 
(Lemma \ref{dual}).
In particular, 
we can find an essential simple closed geodesic $\gamma$ in $\tilde\Sigma\setminus\bo\tilde\Sigma$.

 For (2),  let $g$ be an element of $Q$ that does not lie in an extended boundary subgroup,
and represent the image of $g$ in $\pi_1(\Sigma)$ by an immersed geodesic $\delta$
not parallel to the boundary.

If $\Sigma$ is an orientable  surface, the existence of an  essential simple closed geodesic   implies that $\Sigma$ is filled by such geodesics. In particular, $\delta$ intersects such a simple geodesic $\gamma$,
and $g$ is hyperbolic in the splitting of $G$ dual to $\gamma$. 
Our assumptions guarantee that this splitting  is over $\cala$, so $g$ is not universally elliptic.

If $\Sigma$ is an orbifold, the argument works the same, but we need to check
that simple geodesics fill the orbifold.
Starting with a simple geodesic $\gamma$ as above, 
Lemma 5.3 of \cite{Gui_reading} ensures that there exists another simple geodesic $\gamma_1$  
intersecting $\gamma$ non-trivially. Let $\Gamma_1\subset \Sigma$ be the suborbifold with geodesic boundary filled by $\gamma ,\gamma_1$,
  constructed by projecting the (necessarily $\Lambda$-invariant) subsurface with geodesic boundary filled by the preimage of $\gamma \cup\gamma_1$ in the surface $\Sigma_0$. 
 Applying again Lemma 5.3 of \cite{Gui_reading} to any  component of  $\partial \Sigma\setminus\partial \Gamma_1$, we get a 
larger suborbifold $\Gamma_2$, and the repetition of this process has to stop with some $\Gamma_i=\Sigma$,
which concludes the proof.
\end{proof}

\subsection{When  flexible groups are QH} 

The following theorem collects results from \cite{RiSe_JSJ, DuSa_JSJ, FuPa_JSJ}. 

\begin{thm}\label{thm_lesgens}
Let $G$ be finitely presented.  Let $Q$ be a non-slender 
flexible vertex stabilizer of a JSJ tree
$T_\cala$
 over
$\cala$. In each of the following cases, $Q$ is a QH-subgroup:
\begin{enumerate}
\item  
$\cala$ consists of all finite or  
cyclic
subgroups of $G$.  In this case, $F$ is trivial, and    the underlying orbifold has no mirror. 

\item  
$\cala$ consists of all finite and 2-ended   
subgroups of $G$. In this case, $F$ is    finite.  

More generally, $\cala$ consists of all $VPC_{\le n}$ subgroups of $G$, for some $n\ge1$, 
and $G$ does not split over a $VPC_{\le n-1}$ subgroup. In this case, $F$ is $VPC_{n-1}$.

\item  
$\cala$ consists of all slender  
subgroups of $G$. In this case, $F$ is    slender. \qed
\end{enumerate}
\end{thm}

 These results are  usually stated under the assumption that $G$ is one-ended (in this case, one does not need to include finite subgroups in $\cala$),  
but    they   hold in general (\cite[p.\ 107]{RiSe_JSJ}, see Corollary \ref{cor_unbout}).

\begin{rem}
  In \cite{GL3b}, we will prove that  similar results hold for splittings of   relatively hyperbolic 
groups with small parabolic subgroups over  
  small subgroups.
In \cite{DG2}, one studies the flexible subgroups of the JSJ deformation space of a one-ended hyperbolic group over
the class $\calz$ of virtually cyclic groups with infinite center (and their subgroups).
 Flexible vertices then are QH-subgroups whose underlying orbifold has no mirror. In particular, the fundamental group
of an orbifold with mirrors usually has a non-trivial JSJ deformation space  over $\calz$.
One also studies a (variation of) JSJ splitting over maximal subgroups in $\calz$.
Its flexible groups are \emph{orbisockets}, \ie QH-subgroups without mirrors,
amalgamed to virtually cyclic groups over  maximal extended boundary subgroups (\cite{DG2}, see also \cite{RiSe_JSJ}). 
\end{rem}

Propositions \ref{prop_fue} and \ref{prop_relper} immediately imply:

\begin{cor} \label{corlesgens}
  Under the  assumptions of Theorem \ref{thm_lesgens}, any flexible subgroup of a JSJ decomposition of
$G$ is QH with universally elliptic slender fiber, and all its boundary components are used.
Moreover the underlying orbifold contains an essential simple closed geodesic. \qed
\end{cor}

 For instance, the underlying orbifold cannot be  a pair of pants.

\subsection{Quadratically hanging subgroups are elliptic in the JSJ}
\label{sec_QH_ell}

The goal of this subsection is to prove that,  
under suitable hypotheses, a QH vertex group   of any splitting is elliptic in the JSJ deformation space 
(if we do not assume existence of the JSJ deformation space, we obtain ellipticity in every universally elliptic tree). 

We start with the following fact 
(see \cite[Remark 2.3]{FuPa_JSJ}).

\begin{lem} \label{lem_elhyp} 
Suppose that $T_1$ is elliptic with respect to  $T_2$, but 
 $T_2$ is not elliptic with respect to $T_1$. 
 
Then $G$ splits over a group  which has infinite index in an edge stabilizer of $T_2$.
\end{lem}

\begin{proof} Let $\hat T_1$ be as in Lemma \ref{lem_refinement}.  Let $G_e$ be an edge stabilizer of $T_2$ which is
not elliptic in $T_1$. By Assertion  (3) of Lemma 
\ref{lem_refinement}, it contains an edge stabilizer $J$ of  
$\hat T_1$.  Since $J$ is elliptic in $T_1$ and $G_e$ is not, the index of
$J$ in $G_e$ is infinite.
\end{proof}

\begin{cor} \label{pass}  If $G$ splits over a group $K\in\cala$, but does not split over any infinite index subgroup of $K$, 
then $K$ is elliptic in the JSJ deformation space over $\cala$.  
\end{cor} 

\begin{proof} Take $T_1$ to be a JSJ tree, and $T_2$ a one-edge splitting over $K$. 
\end{proof}

\begin{rem} \label{uni} If $K$ is not universally elliptic, it fixes a unique point in any JSJ tree. Also note that being elliptic or universally elliptic is a commensurability invariant, so the same conclusions hold for groups commensurable to $K$. 
\end{rem}

In \cite{Sela_structure,RiSe_JSJ,DuSa_JSJ}, it is proved that,  if  $Q$ 
 is a  QH-subgroup in some splitting (in the class considered), then $Q$ is 
elliptic in the JSJ deformation space.

This is not true in general, even if $\cala$ is the class of cyclic groups: $F_n$ contains many quadratically hanging subgroups, none of them elliptic in the JSJ deformation space.
This happens because $G=F_n$ splits over groups in $\cala$ having infinite index in each other, 
something which is prohibited by the hypotheses of the papers mentioned above (in \cite{FuPa_JSJ},  $G$ is allowed to split over a subgroup of infinite index in a group in $\cala$, but $Q$ has to be the enclosing group of \emph{minimal splittings}, see  
Definition  4.5 and Theorem  5.13(3) in \cite{FuPa_JSJ}).

A different counterexample will be given in Example \ref{exq} of Subsection \ref{ab}, with $\cala$ the family of abelian groups. 
 In that example the QH-subgroup   
 only has one abelian splitting, which is universally elliptic, so it is not elliptic in the JSJ  space.  

These examples explain the hypotheses in the following general result.

\begin{prop}\label{qhe}
Consider an    $\cala$-tree $T$. Let  $Q $  be a  QH vertex stabilizer of $T$, with fiber $F$ and underlying orbifold $\Sigma$. 
Assume that, if  $\hat J\inc Q$ is the preimage of a 2-ended subgroup   $J\inc \pi_1(\Sigma)$,  
 then $\hat J$  belongs to $\cala$ 
and $G$ does not split over a subgroup of infinite index of $\hat J$.
If $\partial \Sigma=\es$, assume additionnally that   $F$ is universally elliptic.

 Then $Q$ is elliptic in the JSJ deformation space over $\cala$.  
\end{prop}

\begin{rem} \label{rem_qhe}
  We do not assume that the fiber $F$ is slender. 
\end{rem}

\begin{proof}
 We first claim that every boundary component of $ \Sigma$ is used  by an edge of $T$ (see Definition \ref{dfn_qh}). 
Otherwise, we argue as in the proof of Proposition \ref{prop_fue}. Using Lemma \ref{arc}, we construct a splitting of $Q$ relative to its 
incident edge groups over a group $F'$ containing $F$ with index 1 or 2. 
By Remark \ref{arcb}, the   group $F'$ is contained in  the preimage $\hat J$ of a 2-ended subgroup $J\inc \pi_1(\Sigma)$. The splitting of $Q$ extends to a splitting of $G$, and this is a contradiction since $F'$ has infinite index in $\hat J$.

   We deduce that $G$ splits non-trivially over the preimage $\hat B$ of every boundary subgroup $B\inc \pi_1(\Sigma)$. Let  $T_J$ be a JSJ tree,  and assume that $Q$ is not elliptic in $T_J$. 
Since $\hat B\in\cala$ by assumption, $\hat B$ is elliptic in   $T_J$ by Corollary \ref{pass}.
 In particular, $F$ is elliptic in $T_J$  (this follows from our additional assumption if $\partial \Sigma=\es$). As in Proposition \ref{prop_fue},
  the action of $Q$ on the fixed point set of $F$ in $T_J$ factors through an action of $\pi_1(\Sigma)$.  Note that every boundary subgroup is elliptic, 
so Lemma \ref{dual}
implies the existence of a simple closed geodesic $\gamma$ in $\tilde\Sigma\setminus\bo\tilde\Sigma$
(as defined in Subsection \ref{orb}).

 Given such a $\gamma$, recall that $H_\gamma$ is the subgroup of $\pi_1(\Sigma) $ consisting of elements which preserve $\gamma$ and each of the   half-spaces bounded by    $ \gamma$. We let $\hat H_\gamma$ be the preimage of $H_\gamma$ in $Q$.  
It  belongs to $\cala$,
and $G$ has a non-trivial one-edge splitting over $\hat H_\gamma$. 
Lemma 5.3 of \cite{Gui_reading} implies that every essential simple closed  geodesic 
crosses another one. In particular, $\hat H_\gamma$ is not universally elliptic.  
By Corollary \ref{pass}, $\hat H_\gamma$  fixes a   point $c_\gamma\in T_J$.  By Remark \ref{uni}, this point is unique.

\begin{lem} 
 Let $\gamma,\gamma'\inc\tilde\Sigma\setminus\bo\tilde\Sigma$ be simple geodesics. If $\gamma$ and $\gamma'$ intersect, then $c_\gamma=c_{\gamma'}$.
 \end{lem}

\begin{proof} Let $T$ be the Bass-Serre tree of the splitting of $G$ determined by $\gamma'$. It contains an edge $e$ with stabilizer $\hat H_{\gamma'}$. Since $\gamma$ and $\gamma'$ intersect, the group $\hat H_\gamma$ acts hyperbolically on $T$, and its minimal subtree $M$ contains $e$. Let $  T_1$ be a refinement of $T_J$ which dominates $T$, and let $  M_1\inc  T_1$ be the minimal subtree of $\hat H_\gamma$. 

The image of $  M_1$ in $T_J$ consists of the single point $c_\gamma$ (because $  T_1$ is a refinement of $T_J$), and its image by any equivariant map $f:  T_1\to T$ contains $M$. 
Let $  e_1$ be an edge of $  M_1$ such that $f(  e_1)$ contains $e$. The stabilizer $G_{  e_1}$ of $  e_1$ is contained in $\hat H_{\gamma'}$ \emph{with finite index}, so $c_{\gamma'}$ is its unique fixed point. But $G_{  e_1}$ fixes $c_\gamma$, so $c_\gamma=c_{\gamma'}$.
\end{proof}

We can now conclude. First suppose that $\Sigma$ is a surface. Choose a finite set $\Gamma$ of simple closed geodesics which fill $\Sigma$. 
The family $\tilde \Gamma$ of lifts of elements of $\Gamma$ is then connected. 
  Given $\gamma,\gamma'\in\tilde\Gamma$, we can find simple geodesics $\gamma=\gamma_0,\gamma_1,\dots,\gamma_p= \gamma'$ such that $\gamma_i$ and $\gamma_{i+1}$ belong to $\tilde\Gamma$ and intersect.  The lemma implies that $c_\gamma$ is independent of $\gamma\in\tilde\Gamma$. It is therefore 
    fixed by $Q$,  a contradiction. 
    
   The proof is similar when $\Sigma$ is an orbifold, taking $\Gamma$ to be a set of simple geodesics which fill $\Sigma$ as in the proof of Proposition \ref{prop_relper}.   
     \end{proof}

\begin{rem}
Let $Q$ be  a QH vertex  stabilizer as in Proposition \ref{qhe}. Assume moreover that all groups in $\cala$ are slender,
and that $G$ does not split over a subgroup of $Q$ whose
 image in $\pi_1(\Sigma)$ is finite.
Then $Q$ fixes a unique point $v\in T_J$,  so $Q\subset G_v$.
We claim that, \emph{if the stabilizer $G_v$ (hence also $Q$) is universally elliptic, then $v$ is a QH-vertex of $T_J$.}
This is  used in \cite{GL5}.

Let $T$ be an $\cala$-tree in which $Q$ is a vertex group $G_w$. 
First note that $G_v$ is elliptic in $T$, so $G_v=Q$.   We have  to show that, if $e\inc T_J$ contains $v$, then $G_e$ is 
 an extended boundary subgroup of $Q$.
Let $\hat T$ be a refinement of $T_J$ which dominates $T$,  as   in Lemma \ref{lem_refinement}.
Let $\hat w$ be the unique  point of $\hat T$ fixed by $Q$, and let $f:\hat T\to T$ be an equivariant map.  Let $\hat e$ be the lift of $e$ to $\hat T$.

If $f(\hat e)\ne\{w\}$,  then $G_e$ fixes an edge of $T$ adjacent to $w$, so is
an extended boundary subgroup of $Q$. Otherwise, consider
 a segment $x\hat w$, with $f(x)\ne w$, which contains $\hat e$. Choose such a segment of minimal length,
and let $\varepsilon=xy\ne \hat e$ be its initial edge.
We have $G_\varepsilon\inc G_y\inc G_w =Q$, and $G_\varepsilon$ fixes an edge of $T$ adjacent to $w$.
Since $G$ does not split over groups mapping to finite groups in $\pi_1(\Sigma)$,
the image of $G_\varepsilon$ in $\pi_1(\Sigma)$ is a finite index subgroup of a boundary
subgroup $B\inc\pi_1(\Sigma)$.
But we also have $G_\varepsilon\inc G_w=G_{\hat w}$, so that $G_\varepsilon\inc G_{\hat e}=G_e$.
Being slender,  the image of $G_e$ in $\pi_1(\Sigma)$ has to be contained in $B$.
\end{rem}

\section{JSJ decompositions with flexible vertices}
\label{sec_flex}

We shall now consider   JSJ trees with flexible vertices. By Lemma \ref{lem_flex}, the stabilizer of such a
vertex has a non-trivial splitting (relative to the incident edge groups), but not over a universally elliptic 
subgroup.

\subsection{Changing edge  groups}\label{ceg}

First we fix two families of subgroups $\cala$ and $\calb$
and we   compare the   associated JSJ splittings.   For  example:  

$\bullet$ $\cala$ consists of the trivial group, or   the finite subgroups of $G$.
This 
 allows us to reduce to   the case when $G$ is   one-ended  
when discussing flexible
subgroups  (see Corollary \ref {cor_unbout}).

$\bullet$ $\cala$ consists of the finitely generated abelian subgroups of $G$, and $\calb$ consists of
the slender subgroups. This will be useful to describe the abelian JSJ in
Subsection \ref{ab}.

$\bullet$  $G$ is relatively hyperbolic, $\cala $ is the family of parabolic
subgroups, and $\calb$ is the family of elementary subgroups. 

There are now two notions of universal ellipticity, so we  shall distinguish between $\cala$-universal ellipticity (being elliptic in all $\cala$-trees) and $\calb$-universal ellipticity. 

Two trees are \emph{compatible} if they have a common refinement. 

\begin{prop}\label{prop_embo}
 Assume   $\cala\subset\calb$. Let $T_\calb$ be a JSJ tree over $\calb$.
  If there is a JSJ tree over $\cala$, there is one  which is compatible with $T_\calb$. It may be
  obtained    by refining $T_\calb$, and then collapsing all edges whose
stabilizer is not in $\cala$.
\end{prop}

\begin{proof} Let $T_2$ be a JSJ tree over $\cala$. Apply
Lemma \ref{lem_refinement} with $T_1=T_\calb$. Consider an edge $e $ of
$\hat T_1$ whose stabilizer is not in $\cala$. Then $G_e$ fixes a unique
point of $T_2$,   so any equivariant map   from $\hat
T_1$ to
$T_2$ 
is constant on $e$. It  follows that the tree obtained from $\hat T_1$ by
collapsing all edges whose stabilizer is not in $\cala$  dominates
$T_2$. Being $\cala$-universally
elliptic by Assertion  (2) of Lemma \ref{lem_refinement}, it  is a JSJ tree over $\cala$. 
\end{proof}

\begin{prop}\label{prop_embo2}
Assume  that  $\cala\subset\calb$, and 
every $\cala$-universally elliptic $\cala$-tree is $\calb$-universally
elliptic.
\begin{enumerate}
 \item  Let $T_\calb$ be a JSJ tree over $\calb$,  and  let
$\pi :T_\calb\to T _\cala$ be  the map that
collapses all edges whose stabilizer is not in $\cala$. Then:
\begin{enumerate}
\item $T  _\cala$ is a JSJ tree over $\cala$;
\item  if $v$ is a vertex of $T  _\cala$, with stabilizer $G_v$, then 
$\pi^{-1}(v)$  contains a JSJ tree  $T_v$  of
$G_v$ over $\calb_{|G_v}$ relative to the incident edge groups.  
\end{enumerate}
\item Conversely, suppose that $T _\cala$ is a JSJ tree over $\cala$, and for every  stabilizer $G_v$ there exists a JSJ tree  $T_v$  over $\calb_{|G_v}$ relative to the incident edge groups.  Then one can refine  $T _\cala$ using these  trees so as to get  a JSJ tree over $\calb$.
\end{enumerate}
\end{prop}

\begin{proof}  
For the first assertion, let $T$ be an $\cala$-universally elliptic $\cala$-tree. It is  $\calb$-universally
elliptic, so is dominated by $T_\calb$. 
As in the previous proof, one shows that any equivariant map from $T_\calb$ to $T$ factors through $T_\cala$, so $T_\cala$ is a JSJ tree over $\cala$. 
Statement (b) follows from Assertion  (3) of Lemma \ref{lem_flex} (applied over $\calb$).

For the second assertion, we view $T_\cala$ as a $\calb$-universally elliptic  tree.
The proposition then follows from Lemma \ref{lem_JSJrel} (applied over $\calb$).
\end{proof}

\begin{rem}\label{rem_jrel}
This proposition   remains  true in a relative setting,  provided that   one enlarges $\calp_v$ 
  as in Remark \ref{rem_flexrel}.
\end{rem}

The hypothesis of the proposition is satisfied if all
groups in
$\cala
$ are finite (or, more generally, have Serre's property FA). In this case  $T_v$ is simply a (non-relative) JSJ tree over $\calb_{|G_v}$.  In particular, letting $\cala $ consist of the trivial group (resp.\ all finite subgroups), we deduce that one may assume one-endedness of $G$ when studying flexible subgroups.

\begin{cor}\label{cor_unbout}  Let $\calb$ be arbitrary.
 
\begin{enumerate}
\item $G$ has a JSJ deformation space over $\calb$ if and only if each of its non-cyclic free factors does. If so,
every  flexible subgroup of  $G$ is a flexible
subgroup of a non-cyclic  free factor.
\item If $\calb$ contains all finite subgroups, then $G$ has a JSJ deformation space over $\calb$ if and only if $G$  is accessible (see Subsection \ref{Dunw}) and every maximal one-ended subgroup has a JSJ space. If so, every flexible
subgroup of 
$G$    is a flexible subgroup   of a maximal 
one-ended subgroup.
\qed
\end{enumerate}
\end{cor}

\subsection{Peripheral structure of quadratically hanging  vertices}\label{periph}

In this subsection, we study the incident edge groups of a QH vertex stabilizer, and how they depend on the particular JSJ 
tree chosen.

Let  $Q$ be a  QH vertex stabilizer of a JSJ tree $T_\cala$ as in Theorem \ref{thm_lesgens}.
If $T'$ is another JSJ tree, it has a vertex stabilizer equal to $Q$, but possibly
with different incident edge groups. 
For instance, 
it is always possible to modify
$T_\cala$  within its deformation space so that  each incident edge group is 
a maximal  extended boundary subgroup of $Q$  (see Definition \ref{dfn_qh}).
But one loses information in the process (see Example \ref{perd} below).  
 
The relevant structure, which does not depend on the
choice of a JSJ tree,  is the peripheral structure of $Q$, as defined in Section 4 of \cite{GL2}. This is a
finite family $\calm_0$ of conjugacy classes of extended boundary subgroups of $Q$.  In the case at hand,
one may define
$\calm_0$ as follows (see \cite{GL2} for details): a subgroup of
$Q$ represents an element of $\calm_0$ if and only if it fixes an edge in every
JSJ tree, and is maximal for this property.  

We have seen  (Corollary \ref{corlesgens}) 
 that $\calm_0$ uses   every boundary component of $\Sigma $.
Apart from that, the peripheral structure of $Q$ may be fairly arbitrary.
We shall now give examples. In particular, \emph{even when $G$ is one-ended, it is possible for an incident edge
group to meet
$F$ trivially, or (in the slender case) to have trivial image in $\pi _1(\Sigma )$.}

\paragraph{Construction.}
 Let  $Q$ be an extension $1\to F\to Q \to \pi
_1(\Sigma )\to 1$  with $F$   slender and $\Sigma $ a compact orientable surface  
(with genus $\ge 1$, or with at least 4  boundary components). Let $H_1,\dots,H_k$ be a finite
family of infinite   extended boundary subgroups of $Q$  as defined in Definition \ref{dfn_qh} (note that they are slender).   
For each boundary subgroup $B$ of
$\pi _1(\Sigma )$, there should be an $i$ such that a conjugate of $H_i$ maps onto
a finite index subgroup of $B$ (i.e.\  every boundary component is used  in the sense of  Definition \ref{dfn_qh}). 
Let $R_i$ be a non-slender finitely presented group with Serre's property FA, for
instance $SL(3,\Z)$. 
We define a finitely presented group $G$ by amalgamating $Q$ with   $K_i=H_i\times R_i$ over $H_i$ for
each $i$; in other words, $G=((Q*_{H_1}K_1)*\dots)*_{H_k} K_k$.

\begin{lem}\label{lem_amal}
The Bass-Serre tree $T$ of the amalgam defining $G$ is a slender JSJ tree,  
$Q$ is a flexible subgroup, and $G$ is one-ended.  If no $H_i$ is conjugate  in $Q$ to a subgroup of
  $H_j$ for $i\ne j$, the peripheral structure $\calm_0$ consists of
the conjugacy classes of the $H_i$'s.
\end{lem}

\begin{proof}  We work over the family $\cala$ consisting of all slender subgroups. 
Let
$T'$ be any tree with an action of $G$. Each
$R_i$ fixes a unique point, and this point is also fixed by $H_i$. In particular,
$H_i\times R_i$, $H_i$, and $T$, are universally elliptic. 
To prove that $T$ is a JSJ tree, it suffices to see that $Q$ is elliptic in any universally elliptic $\cala$-tree $T'$.

By Proposition \ref{prop_fue} (\ref{it_fue}), $F$  is  
universally elliptic.  If $Q$ is not elliptic in $T'$, then by   Proposition \ref{prop_fue} (\ref{it_factor})
the action of $Q$ on its minimal subtree    $T_Q\subset T'$
   factors through a nontrivial action of $\pi _1(\Sigma )$ with slender (hence cyclic) edge stabilizers. Since $H_i$, hence every
boundary subgroup of $\pi _1(\Sigma )$, is elliptic, the action  is dual to a
system of disjoint  essential simple closed geodesics on $\Sigma$  by Lemma \ref{dual}.
  By Proposition \ref{prop_relper}, 
no edge stabilizer of $T_Q$ is universally elliptic,  contradicting  universal ellipticity of $T'$.

Thus, $T$ is a JSJ tree, and 
 $Q$ is flexible because $\Sigma $ was chosen to contain
intersecting simple closed curves.

By Proposition  \ref{prop_embo2}, one obtains a JSJ tree of $G$ over finite groups
by collapsing all edges of $T$ with infinite stabilizer. Since each $H_i$ is infinite,
this JSJ is trivial, so $G $ is one-ended. 

The assertion about $\calm_0$ follows from the definition of $\calm_0$ given in
\cite{GL2}. 
\end{proof}

\begin{example}\label{perd}
Let $\Sigma $ be a punctured torus, with fundamental group  $\langle a,b\rangle$.
Write $u= [a,b]$. Let $Q=F\times\langle a,b\rangle$, with $F$ finite and
non-trivial. Let $H_1=\langle F, u^2\rangle$ and $H_2=\langle u\rangle$. The tree
$T$ is a JSJ tree over 2-ended  (and finite) groups. The peripheral structure of $Q$ consists of
two elements, though $\Sigma $ only has one boundary component. There is a JSJ
tree $T'$ such that incident edge groups are conjugate to $\langle F, u \rangle$
(the quotient $T'/G$ is a tripod), but it does not display the peripheral structure of $Q$.
\end{example}

\begin{example} Let $\Sigma , a,b, u$ be as above. Again write $Q=F\times\langle
a,b\rangle$, but now $F=\langle t\rangle$ is infinite cyclic. Let $H_1=\langle  
u\rangle$ and $H_2=\langle t\rangle$. Then
$H_1$ meets
$F$ trivially, while $H_2$ maps trivially into $\pi _1(\Sigma )$.

\end{example}

\subsection{Flexible vertices of  abelian JSJ decompositions} 
\label{ab}
  
We have described flexible subgroups over cyclic groups, 2-ended
groups, slender groups  (see Theorem \ref{thm_lesgens}). Things are more complicated over abelian groups. 

The
basic reason  is the following: if a group $Q$ is an extension $1\to F\to Q \to \pi
_1(\Sigma )\to 1$ with $F$  a finitely generated abelian group and $\Sigma $ a
surface, a    splitting of
$\pi _1(\Sigma )$ along a simple closed curve induces a splitting of $Q$ over a   subgroup which is
slender (indeed polycyclic) but not necessarily abelian. In the language of \cite{FuPa_JSJ}, the enclosing graph decomposition of two splittings over abelian groups is not necessarily over abelian groups. 

In fact, we shall now   construct examples showing:   
\begin{prop}
\begin{enumerate}
\item  {Flexible subgroups of abelian JSJ
trees are not always  slender-by-orbifold groups.} 
\item {One cannot
always   obtain
an abelian JSJ tree by collapsing edges in a slender JSJ tree.}
\end{enumerate}
\end{prop}

By Proposition \ref{prop_embo}, 
 one can obtain an abelian JSJ tree by refining and collapsing  a slender JSJ tree. The point here is that collapsing alone is not always sufficient. 
 We will see, however, that things
change if
$G$ is assumed to be  CSA  (see Proposition \ref{prop_JSJ_CSA}).

We use the same construction as in the previous subsection, but now   $\pi _1(\Sigma )$ will act non-trivially on the fiber $F$.

\begin{example}  In this example $F\simeq\Z$. Let $\Sigma $ be obtained by
gluing a once-punctured torus to one of the boundary components of a pair of
pants. Let
$M$ be a circle bundle over $\Sigma $ which is trivial over the punctured torus
but non-trivial over the two boundary components of $\Sigma $. Let $Q=\pi
_1(M)$, and
$H_1,H_2$ be the fundamental groups of the components of $\partial M$
(homeomorphic to Klein bottles). Note that $H_1,H_2$ are non-abelian. Construct
$G$ by amalgamation  with $H_i\times R_i$ as above. 
We claim that \emph{the abelian JSJ
decomposition of
$G$ is trivial, and $G$ is flexible (but not slender-by-orbifold). }

We argue as in the proof of Lemma \ref{lem_amal}.
We know that $H_1\times R_1$ and $H_2\times R_2$ (hence also $F$) are universally elliptic.  If  $T$ is any tree with abelian edge stabilizers, the  action of $Q$  on its minimal subtree factors through $\pi_1(\Sigma)$, and the action of $\pi_1(\Sigma)$ is dual to a system of simple closed curves
 (Proposition \ref{prop_fue} (\ref{it_factor}) and Lemma \ref{dual}). 
But not all curves give rise to an abelian splitting: they have to be ``positive'', in the sense that  the bundle is trivial over them.

To
prove that none of these splittings is universally elliptic, hence that the abelian JSJ space of $G$  is trivial,   it suffices to see that any positive
curve intersects (in an essential way) some other positive curve. This is true for
the curve $\delta $ separating the pair of pants from the punctured torus (one easily
constructs a positive curve meeting $\delta $ in 4 points). It is  also 
true for
curves meeting $\delta $. Curves disjoint from $\delta $ are contained in the
punctured torus, and the result is true for them.

It is clear that $G$ is flexible. If one performs the construction adding a third group $H_3=F$, then $G$
becomes a flexible vertex group in a group $\hat G$ with non-trivial abelian
JSJ. 
\end{example}

\begin{example}  \label{exq}
Now $F=\Z^2$. Let $\Sigma $ be a surface of genus $\ge2$ with two boundary
components $C_1,C_2$. Let  $\gamma $ be  a simple closed curve 
separating $C_1$ and $C_2$. Let $\Sigma '$ be the space obtained from $\Sigma
$ by collapsing $\gamma $ to  a point. 
Map $\pi _1(\Sigma )$ to  $SL(2,\Z)\inc\Aut(\Z^2)$ by projecting to  $\pi
_1(\Sigma ')$ and embedding the free group  $\pi _1(\Sigma ')$ into $SL(2,\Z)$.
Let $Q$ be the associated semi-direct product $\Z^2\rtimes \pi _1(\Sigma )$, and
$H_i=\Z^2\rtimes \pi _1(C_i)$.  Construct $G$ as before. 

Abelian splittings of $G$ now
come from simple closed curves on $\Sigma $ belonging to the kernel of $\rho: \pi _1(\Sigma )\to SL(2,\Z)$. 
But it is easy to see that $\gamma $ is the only
such curve. It follows that  the one-edge splitting dual to
$\gamma $ is an abelian JSJ  decomposition of $G$. It has two rigid vertex groups. It cannot be obtained by collapsing a
slender JSJ splitting. 
\end{example}

\begin{rem}
In order to obtain an abelian JSJ tree from a slender
JSJ tree, one must in general      refine the tree and then
collapse edges with non-abelian stabilizer (Proposition \ref{prop_embo}). There is some
control over how a flexible vertex group $Q$  is   refined. Suppose for instance that
$Q=\Z^n\rtimes \pi _1(\Sigma )$ with $\Sigma $ an orientable  surface. As in the previous
example, the refinement uses curves $\gamma  $ in the kernel of $\rho: \pi _1
(\Sigma )\to\Aut(\Z^n)$ which do not intersect other curves in the kernel. One
may check that there is at most one such curve, and
it is separating. \end{rem}

We now show that the abelian JSJ is very easy to describe when $G$ is CSA.  
Recall  that $G$ is CSA if the centralizer of any non-trivial element is abelian and  malnormal.

\begin{prop}\label{prop_JSJ_CSA} Let $G$ be  a   torsion-free, finitely
presented CSA group. 
\begin{enumerate}
\item
One obtains an abelian JSJ tree $T $ by collapsing all edges
with non-abelian stabilizer in a slender JSJ tree. 
\item All non-abelian flexible subgroups of $T $ are   QH-subgroups with trivial fiber:
they are isomorphic to $\pi _1(\Sigma )$, with
$\Sigma $ a surface and all incident edge groups  contained in boundary subgroups.  
\end{enumerate}
\end{prop}

See \cite{GL3b} 
for abelian JSJ decompositions of \emph{finitely generated} torsion-free CSA groups.

\begin{proof}   Denote by $\cala$ and $\cals$ the families of abelian and slender subgroups of $G$ respectively.    If $G$ has infinitely generated abelian subgroups, one does not have $\cala\inc\cals$. But 
any $\cals$-universally
elliptic subgroup is  $\cala$-universally
elliptic by Corollary \ref{cor_fg}.

Let $T_\cala$, $T_\cals$ be JSJ trees. 
The key point is to show that $T_\cals$ dominates 
$T_\cala$. This implies   Assertion   (1)   as in the proof of Proposition \ref{prop_embo}:   any map  $T_\cals\to T_\cala$ sends edges with non-abelian stabilizer  to points,  so collapsing these edges yields an abelian JSJ tree.

To prove  that $T_\cals$ dominates 
$T_\cala$, we consider 
 a vertex stabilizer $Q=G_v$ of $T_\cals$ which is not $\cala$-universally
elliptic (hence not $\cals$-universally
elliptic), and we show that it is elliptic in $T_\cala$. 

We first show that $Q$ is either an abelian group or a surface group (a QH-subgroup with $F$ trivial). 
There are two cases. 
If $Q$   is slender, let $S$ be an $\cala$-tree on which $Q$ acts non-trivially.  
 By slenderness,  $S$ contains
an $Q$-invariant line, and    {$Q$ is abelian} by the CSA property as it maps onto  $\Z$ or an infinite  dihedral group
with  abelian kernel. 
If $Q$ is not slender, it is a QH-subgroup  with slender fiber $F$ by  Theorem   \ref{thm_lesgens}. 
We prove that $F$ is abelian, hence trivial by the CSA property. We have seen that $F$ is $\cals$-universally
elliptic (Proposition \ref{prop_fue}). If $F$ were not abelian, it would fix a unique point in every $\cala$-tree. 
This point would be fixed by $Q$ because $F$ is normal, contradicting the fact that  $Q$ is not $\cala$-universally elliptic.

 We  can now show that $Q=G_v$ is elliptic in $T_\cala$.   Since  
  all its slender subgroups are abelian,  the JSJ deformation space of $Q$ relative to the incident edge groups  is the same over $\cala$ as over $\cals$.
 Note that these edge groups are abelian. Applying   Assertion  (3) of Lemma \ref{lem_flex} 
 with $T$ the tree   obtained from $T_\cals$ by collapsing all edges with non-abelian stabilizer, we see that $Q$ is elliptic in the JSJ space over $\cala$ since it is elliptic over $\cals$.  This proves  the first assertion of the proposition.

We now prove  the second assertion. By 
(1), we can assume that $T_\cala$ is obtained from $T_\cals$ by collapsing. 
No edge adjacent to   
a QH-vertex of $T_\cals$ with trivial fiber is collapsed, so such a
vertex remains a flexible vertex stabilizer  of  $T_\cala$, with the same  incident edge groups. 
 We have seen that   all other vertex stabilizers of $T_\cals$ are abelian or $\cala$-universally elliptic.
Thus  a vertex stabilizer  of
$T_\cala$ which is not an  abelian group   or a surface group is    the fundamental
group of a graph of groups with non-abelian edge groups and $\cala$-universally elliptic vertex
groups. Such a group is $\cala$-universally elliptic. It follows that all non-abelian flexible subgroups of $T_\cala$ are surface groups. 
\end{proof}

\subsection{Relative JSJ decompositions}\label{relb}

  Fix   a finitely
presented group $G$ and a family $\cala$. In Subsection \ref{rela} we have shown the existence of  the JSJ deformation space of  $G$  relative to a finite set  $\calh=\{H_1,\dots,H_p\}$  of
finitely generated subgroups.  
We now give an alternative construction, using   the (absolute) JSJ
space of another group $\hat G$ obtained by the filling construction of  Subsection \ref{periph}.  This will  allow us in the next subsection  to 
extend  the description of flexible vertices of Theorem
\ref{thm_lesgens}   to the relative case.

Let $G$  and $\calh$ be as above. 
As  in Subsection \ref{periph}, 
we define a group
$\hat G$ by amalgamating   with  $K_i=H_i\times R_i$ over $H_i$, where $R_i$ is a 
finitely presented group with property FA. The group $\hat G$ is the fundamental group of a graph of groups with one central vertex $v$ (with vertex group $G$), and edges $e_i=vv_i$ with $G_{e_i}=H_i$ and
$G_{v_i}=H_i\times R_i$. It  is finitely presented. We denote by $T$ the Bass-Serre tree of this amalgam. 

\begin{rem} 
This construction   
can be adapted  to  the more general setting where 
$G$ is finitely presented relative to $\{H_1,\dots,H_p\}$, and each $H_i$ is finitely generated and recursively presented.
  Since $H_i$ embeds   into a finitely presented group $H'_i$, 
 one can take for $\hat G$ the amalgam of $G$ with $H'_i\times R_i$ over $H_i$, a finitely presented group.
\end{rem}

Fix a family $\calb$ of
subgroups of $\hat G$ such that $\calb_{|G}=\cala$ and $R_i\notin \calb$,  for instance 
the family of subgroups of $\Hat G$ having a  conjugate in $\cala$ (note that  two subgroups of $G$ which are conjugate in $\hat G$ are also conjugate in $G$  because $H_i$ is central in $K_i$,  so this family $\calb$ is stable under conjugation).

We also define a family $\calb_\calh$,  
by adding to $\calb$  
 all subgroups of $\Hat G$ 
 having a conjugate contained in some
$H_i$.  Note that  
 $R_i\notin \calb_\calh$, since  dividing $\hat G$ by the normal closure of $G$
kills every $H_i$ but does not kill $R_i$. 

\begin{lem} \label{lem_relu} $H_i\times R_i$ is $\calb_\calh$-universally elliptic. 
A subgroup $J\subset G$ is $(\cala, \calh)$-universally elliptic if and only if $J$ (viewed as a subgroup of $\hat G$) is
$\calb $-universally elliptic.
\end{lem}

\begin{proof} Consider any $\calb_\calh$-tree. The group $R_i$ fixes a point, which is
unique since $R_i\notin \calb_\calh$. This point is also fixed by $H_i$ since $H_i$
commutes with $R_i$.   Since $T$ is $\calb $-universally elliptic, the second assertion follows from Lemma \ref{lem_flex}.
\end{proof}

Given an $(\cala,\calh)$-tree $T_v$, with an action of $G$, 
we may use $T$ to extend the action of $G$ on $T_v$ to an action of $\hat G$ on a
 tree $\tilde T$ containing $T_v$,  as in Remark   \ref{rem_bas} 
(at the graph of groups level, this amounts to gluing edges amalgamating $H_i$ to $H_i\times R_i$ onto $T_v/G$).
To get a $\calb$-tree $\hat T$, we collapse the orbit of the edge stabilized by $H_i$ whenever $H_i\notin \cala$. 
We say that $\hat T$ is an \emph{extension} of $T_v$ to $\hat G$.

 \begin{lem} \label{lem_some}
 Some JSJ tree $T_\calb$ of $\Hat G$ over $\calb$   is an extension 
of a JSJ tree $T_\cala$ of $G$ over $\cala$ relative to $\calh$.
 \end{lem}

 \begin{proof}  Being finitely presented, $\hat G$ has a JSJ deformation space over $\calb_\calh$ by Theorem \ref{thm_exist_mou}.
The Bass-Serre tree $T$  is
$\calb_\calh$-universally elliptic, so by Lemma \ref{lem_rafin} it may be refined  to a
JSJ tree of $\hat G$ over
$\calb_\calh$. By Proposition \ref{prop_embo}, this JSJ tree is compatible with some
JSJ tree
$T_\calb$ of $\hat G$ over $\calb$.

One obtains $T_\calb$ from $T$ by refining, and then collapsing all edges with stabilizer not in $\calb$. By universal ellipticity of   $H_i\times R_i$, no refinement takes place over the vertex $v_i$, so only $v$ is refined. The edge $e_i$ is collapsed if and only if $H_i\notin\cala$. This shows that   $T_\calb$   is an extension of an  $(\cala,
\calh)$-tree
$T_\cala$. 

The tree
$T_\cala$ is $(\cala, \calh)$-universally elliptic by Lemma \ref{lem_relu}. If it is not maximal, it may be
refined as an
$(\cala,
\calh)$-universally elliptic tree. An  extension of this refinement is $\calb$-universally
elliptic by Lemma \ref{lem_relu}, and contradicts the maximality of $T_\calb$ as a
$\calb$-universally elliptic tree.  This shows that $T_\cala$ is a JSJ tree of $G$ relative to $\calh$.
\end{proof}

This gives another proof of Theorem \ref{thm_exist_mou_rel}:
 \begin{cor} The JSJ deformation space 
 of $G$ over $\cala$ relative to $\calh$ exists.  \qed
 \end{cor}

\subsection{Relative QH-subgroups}

In a relative setting, one needs to slightly modify the definition of a QH-subgroup by taking groups of $\calh$ into account.

\begin{dfn}\label{dfn_relqh}   Given $G$,  $\cala$, and a family of subgroups $\calh$, one says that $Q$ is a \emph{relative QH-subgroup} (over $\cala$, relative to
  $\calh$) if: 
  \begin{enumerate}
  \item $Q$ is a vertex stabilizer of an $(\cala,\calh)$-tree;
  \item $Q$ is an extension $1\to F\to Q \to
  \pi_1(\Sigma )\to 1$, with $\Sigma$ a hyperbolic 2-orbifold; 
  \item each incident edge group, and each intersection of a
  conjugate of a group in $\calh$ with $Q$, is an extended boundary subgroup as in Definition \ref{dfn_qh}.
\end{enumerate}
  A boundary component $C$ of $\Sigma$ is \emph{used}  if    there is an incident edge group, 
or a conjugate of a group in $\calh$, whose image     in $\pi_1(\Sigma)$ is contained in $\pi_1(C)$ with   finite index.
\end{dfn}

\begin{rem} 
As in Remark \ref{rem_flexrel}, one can replace 
the assumption on $\calh$  in (3) 
by the following one: for any $H\in \calh$ and any $g\in G$ such that $H^g\subset Q$,
then $H^g$ is an extended boundary subgroup of $Q$.
\end{rem}

\begin{rem}\label{rem_abs2rel}
  Propositions \ref{prop_fue} and \ref{prop_relper} 
  apply in a relative setting (the meaning of \emph{QH} and of \emph{used} being understood
  as in Definition \ref{dfn_relqh}).
\end{rem}

The description of flexible vertices of JSJ decompositions over classes of slender groups
occuring in  Theorem \ref{thm_lesgens} can be extended to the relative case.

\begin{thm}\label{thm_lesgensrel}
Let $G$ be finitely presented, let $\cala$ be a class of subgroups as below,
and  let $\calh=\{H_1,\dots,H_p\}$ be a finite set of
finitely generated subgroups.   
Let $T$ be a JSJ decomposition of $G$ over $\cala$ relative to $\calh$.

In each of the following cases, the flexible vertices of $T$ with non-slender stabilizer 
are relative QH-subgroups, with fiber $F$ as in Theorem \ref{thm_lesgens},
and every boundary component is used: 
\begin{enumerate}
\item  
$\cala$ consists of all finite or 
cyclic
subgroups of $G$.     

\item  
$\cala$ consists of all finite and 2-ended   
subgroups of $G$. 

More generally, $\cala$ consists of all $VPC_{\le n}$ subgroups of $G$, for some $n\ge1$, no $H_i$ is $VPC_{\le n-1}$, 
and $G$ does not split over a $VPC_{\le n-1}$ subgroup relative to $\calh$.

\item  \label{it_slender}
$\cala$ consists of all slender  
subgroups of $G$. 
\end{enumerate}
\end{thm}

  Let $G$ be hyperbolic relative to a family of finitely generated subgroups $\calh=\{H_1,\dots,H_p\}$,  as in Subsection \ref{parab}. 
Recall that $H<G$ is \emph{elementary} if it is
 finite, virtually cyclic, or contained in a parabolic subgroup.  
In \cite{GL3b}, we will study the elementary JSJ deformation space relative 
  to any family $\calh$  containing every $H_i$ which is not small, and we will show that 
 non-elementary flexible subgroups are relative QH-subgroups with finite fiber.

\begin{proof}[Proof of Theorem \ref{thm_lesgensrel}]  
Choose
$R_i$ so that it does not embed into a slender-by-orbifold group, for instance
$R_i=SL(3,\Z)$. Construct $\hat G$ as above,  and 
let $\calb$ be the natural class of groups extending $\cala$ on $\hat G$ (for instance, the class of slender groups of $\Hat G$ in (\ref{it_slender})). 
 Using Lemma \ref{lem_some}, consider a JSJ tree $T_\cala$ of $G$ relative to $\calh$ as above,
such that some extension $\Hat T_\cala=T_\calb$ is a JSJ tree of $\Hat G$.
Recall that $\Tilde T_\cala/\Hat G$ is obtained by attaching edges to $T_\cala/G$ amalgamating $H_i$ to $H_i\times R_i$,
and that $\Hat T_\cala=T_\calb$ is   obtained by collapsing edges of $\Hat T_\cala$ with stabilizer not in $\calb$ 
(i.e.\ not slender in  case (\ref{it_slender})).

By Theorem \ref{thm_lesgens}, the non-slender flexible vertex stabilizers of $T_\calb$ are QH-subgroups (in the $VPC_{ n}$ case, we must check  that $\hat G$ does not split over a $VPC_{\le n-1}$ subgroup; this holds because in any tree with  $VPC_{\le n-1}$ edge stabilizers  $H_i\times R_i$ is elliptic, so   $G$ is elliptic, and all these groups fix the same point because   no $H_i$ is $VPC_{\le n-1}$).

Let $Q=G_v$ be a non-slender flexible vertex stabilizer of $T_\cala$.
We denote by $\Tilde v$ and $\Hat v$ the  vertex corresponding to $v$ in $\Tilde T_\cala$ and $\Hat T_\cala$.
The stabilizer of $\Tilde v$ is $Q$, and the stabiliser of $\Hat v$ is some $\Hat Q\supset Q$.

First, $\Hat Q$ is  flexible since any splitting of $G$ over $\cala$ relative to $\calh$ extends to a splitting of $\hat G$ over $\calb$. 
By Theorem \ref{thm_lesgens}, $\Hat Q$ is an (absolute) QH vertex group of $T_\calb$.
Since $R_i$ does not embed into a slender-by-orbifold group, no edge of $\Tilde T_\cala$ incident on $\Tilde v$ is collapsed in $\Hat T_\cala$.
It follows that $\Hat Q=Q$.
Since $\Hat v$ is a QH vertex, edge groups incident on $\Hat v$ are extended boundary subgroups of $Q$.
This implies that each conjugate of $H_i$ intersects $Q$ in an extended boundary subgroup of $Q$.
Thus, $Q$ is a relative QH vertex group of $T_\cala$.

The fact that every boundary component is used is a consequence of Proposition \ref{prop_fue} (\ref{it_used}), extended
to the relative case by Remark \ref{rem_abs2rel}.
\end{proof}

\small


\def\cprime{$'$} \newcommand{\noopsort}[1]{}

\begin{flushleft}
Vincent Guirardel\\
Institut de Recherche Math\'ematique de Rennes\\
Universit\'e de Rennes 1 et CNRS (UMR 6625)\\
263 avenue du G\'en\'eral Leclerc, CS 74205\\
F-35042  RENNES Cedex\\
\emph{e-mail:}\texttt{vincent.guirardel@univ-rennes1.fr}\\[8mm]

Gilbert Levitt\\
Laboratoire de Math\'ematiques Nicolas Oresme\\
Universit\'e de Caen et CNRS (UMR 6139)\\
BP 5186\\
F-14032 Caen Cedex\\
France\\
\emph{e-mail:}\texttt{levitt@math.unicaen.fr}\\
\end{flushleft}

 \end{document}